\documentclass[11pt, a4paper]{article}

\usepackage[comma, authoryear, numbers]{natbib}
\usepackage{amsmath,amsfonts,amsthm,bm}
\usepackage[parfill]{parskip}
\usepackage{hyperref}
\usepackage{graphicx}
\usepackage[a4paper, margin=25mm]{geometry}
\usepackage{bm}
\usepackage{bbm}
\usepackage{breqn}
\usepackage[toc,page]{appendix}
\usepackage{tabularx}
\allowdisplaybreaks
\usepackage{multirow}
\usepackage{longtable}
\usepackage{booktabs}
\usepackage{bibentry}
\usepackage[multiple]{footmisc}
\usepackage{algpseudocode}
\usepackage{algorithmicx}
\usepackage{algorithm}
\usepackage{subcaption}
\usepackage{etoolbox}
\usepackage{dsfont}
\usepackage{breqn}
\usepackage{tabulary}
\usepackage{lscape}
\date{\today}

\newtheorem{theorem}{Theorem}
\newtheorem{lemma}{Lemma}
\newtheorem{remark}[theorem]{Remark}
\newtheorem{observation}[theorem]{Observation}

\newcommand{\m}[1]{\mathcal{#1}}

\usepackage{amsthm}
\usepackage[dvipsnames]{xcolor}

\usepackage{tocbibind}
\usepackage{comment}

\usepackage{mathtools}

\title{\textbf{Algorithms for min-buying in networks}}

\date{}

\author{Aaditya Bhardwaj\footnote{STOR-i Centre for Doctoral Training, Lancaster
    University, United Kingdom. Email:
    \href{mailto:a.bhardwaj2@lancaster.ac.uk}{a.bhardwaj2@lancaster.ac.uk}}\ , Ben Black\footnote{STOR-i Centre for Doctoral Training, Lancaster
    University, United Kingdom. Email:
    \href{mailto:benb17995@hotmail.co.uk}{benb17995@hotmail.co.uk}}\
     , Trivikram
    Dokka\footnote{Advanced Analytics Group, Air Products Plc, United Kingdom. Email:
    \href{mailto:trivikram.dokka@yahoo.co.uk}{trivikram.dokka@yahoo.co.uk}}\ , Christopher
    Kirkbride\footnote{Department of Management Science, Lancaster University
    Management School, United Kingdom. Email:
    \href{mailto:c.kirkbride@lancaster.ac.uk}{c.kirkbride@lancaster.ac.uk}}
    \footnote{Corresponding author.}\
     }

\begin{document}
\setcitestyle{nosort}

\maketitle
\begin{abstract}
 The paper is motivated by pricing decisions faced by forecourt fuel retailers across their outlets on a road network. Through our modelling approach we are able adapt the network structure to a bipartite graph with demand nodes representing volumes of fuel from customers using a specific route that connects to the seller's outlet nodes that intersect that route on the network. Customers may have their demand satisfied at the lowest priced competitor on their route. However, the seller can satisfy some or all of this demand by matching or beating this price via one of their outlets intersecting the route. We give a practical extension to min-pricing by considering a binary logit variant for buyers evaluating the choice between two sellers.  We derive two MIP formulations for min-buying in the case of general demand. We also propose several constructive heuristics, based on insertion and selection operations, suitable for problem instances beyond the scope of the exact methods. The performance of models and algorithms are evaluated in a numerical study and develop insights from the results. Importantly, we are able to highlight the value of price-matching decisions under buyer demand sensitivity.
\end{abstract}
\small{\textbf{Keywords:} Pricing, networks, mixed integer programming, heuristics, binary logit choice.}

\normalsize
\pagenumbering{arabic}
\section{Problem Description}
We study the following pricing problem: Consider a bipartite graph $B=(\mathcal{O},\mathcal{N},L) $, where $\mathcal{N}$ are demand nodes and $\mathcal{O}$ are supply or outlet nodes.  Outlet nodes are owned by the seller. Let $\mathcal{O}_e$ be the set of outlets that are connected to demand node $e$, via $L$. Demand originating at node $e$, $d_e$, can be satisfied at a fixed (per-unit) price $c_e$. Alternatively, $d_e$ can be satisfied at outlet $f\in \mathcal{O}_e$ at (per-unit) price $p_f$, if $p_f \leq c_e$. However, the proportion of $d_e$ satisfied by $f$  when $p_f=c_e$ is $D^{PM}_{ef}$ (Price Match); when $p_f<c_e$ is $D^{PW}_{ef}$ (Price War). Without loss of generality, the outlet with lower index will satisfy the demand in case of price ties. Note that we have $D^{PM}_{ef} = D^{PW}_{ef} = 0$ when $p_f>c_e$, $p_f > p_{f'}$ ($f\neq f'$) and when $p_f=p_{f^*}$  where $f^*$ is an outlet of lower index than $f$. Also, $D^{PM}_{ef} \cdot D^{PW}_{ef} = 0$. The objective of the seller is to maximize revenue by setting prices at outlets. That is, to solve the following problem:
\begin{equation*}
   \max_{\textbf{p}} \sum_{e\in \mathcal{N}, f\in \mathcal{O}_e} \left( D^{PM}_{ef} + D^{PW}_{ef} \right) \cdot p_f.
\end{equation*}
We assume that price index set $M$ is discrete with the $m^{th}$ price being $b(m)$. We denote the largest price less than $c_e$ as $c_e^-$.

It is common to study pricing problems with an imposed price order, also referred to as a \textit{price ladder}. An ordering $S$ of outlets $\mathcal{O}$ is said to satisfy the price ladder constraint if prices are non-decreasing within this ordering, that is, $p_{S_1}\leq \ldots \leq p_{S_{|\mathcal{O}|}}$.

\subsection{Motivating real life problem}

We are motivated by the following pricing problem faced by forecourt fuel retailers, which we will refer to as FPP in the following. In FPP, retailers are faced with decisions about the prices to be fixed at each of their outlets. Typically, prices are revised periodically, with periodicity changing from a day to several days. Prices stay fixed for at least a day. The fuel demand at an outlet depends on traffic in the network. That is, demand at a particular outlet $f$ arises from customers travelling between an origin $s$ and destination $t$ if there is route between $s$ and $t$ that passes through $f$. Therefore, every route that intersects $f$ is a demand source. Consequently, a customer taking route $r_f$ can choose between outlet $f$ or any other outlet that intersects $f$, seller owned or otherwise. To connect with our problem definition, every route that generates demand will be grouped as $\mathcal{N}$ and all outlets are grouped as $\mathcal{O}$. The edges in $B$ naturally form as a consequence of outlets on routes. For the sake of simplicity, we omit some details of FPP, such as, 
\begin{itemize}
    \item multiple non-seller outlets: we abstract this case by having one outlet per route which is not owned by seller \textit{aka} competitor owned outlets,
    \item we abstract demand at $e\in \mathcal{N}$ to already incorporate the combination of route choice probability and total demand between the respective origin and destination.
\end{itemize}
Both these choices are made for the sake of simplicity and all aspects of our analysis can be applied when these are incorporated explicitly. 

\subsection{Literature review}

In pricing problems, a retailer seeks to price a collection of items to maximise their revenue or profit. Customers are characterized by a set of items (or bundles of items) they are interested in, and a budget for each item (or bundle) in their set. The customers' decision-making process involves selecting an item from their preference set that lies within their budget. The process typically adheres to one of four models: min-buying, max-buying, max-gain and rank-buying. Min-buying leads the customer to choose the item of lowest price. Conversely, max-buying customers choose the item of maximum price. In the max-gain model customers select the item that maximises their utility. This is defined as the difference between the customer's valuation of the item and the item's price. Finally, in the rank-buying model, customers employ a preference ranking for the items in their set and select the item of maximum rank within their budget.

Inspired by the opportunity to leverage customer preference data from product recommendation engines within price optimization, \cite{PaatRusmevichientong} introduced a non-parametric multi-product pricing model. Here, customers have a fixed (uniform) budget and buying behaviour is either rank- or min-buying. It is shown that the problem is NP-complete and, under a price-ladder constraint that specifies an ordering of the prices of the items, they demonstrate that it can be solved in polynomial time.  \cite{aggarwal2004algorithms} relax the uniform-budget, letting customers identify a maximum price they would pay for each item, and present algorithms for variants of the model. Under the price ladder constraint they provide a polynomial time approximation scheme (PTAS) for the rank-buying model that can be extended to max-buying, and a 4-approximation for max-buying with a limited inventory of items. For max-buying without the price ladder constraint a 1.59-approximation is given and this formulation is noted to be applicable to a max-profit objective in which the seller incurs some cost for the items. For all models they consider, an approximation that is logarithmic in the number of customers is provided. These algorithms are shown to be the best possible by \cite{BuyingCheapisexpensive_Briest}. \cite{Guruswami2005} present max-gain selection under the condition that customers' purchasing is envy-free, i.e. given the pricing no customer would prefer a different bundle of items from that purchased. In the case of unlimited supply, all purchasing is naturally envy-free. Two customer classes are considered, unit demand customers that seek to buy at most one item and single-minded customers who seek to buy a specific bundle and will not purchase anything if the bundle is not affordable. They show that both cases are APX-hard and provide an $\mathcal{O}(log|\mathcal{N}|)$ approximation algorithm for the unit demand with limited supply case and an $\mathcal{O}(log|\mathcal{N}|+|\mathcal{F}|)$ approximation algorithm for the single-minded buyer with unlimited supply. Note that, here, $\mathcal{N}$ denotes the set of customers and $\mathcal{F}$ the set of items. \cite{UniformBudget_Briest} studied the unit-demand min-buying pricing problem where customers have a uniform budget for items. This is a special case of the unit demand envy-free pricing problem since the lowest priced item under consideration will maximize the customer's utility. They showed that the problem cannot be approximated within $\mathcal{O}(log^{\epsilon}|\mathcal{N}|)$ for some $\epsilon > 0$. For the original rank buying model of \cite{PaatRusmevichientong}, the first exact formulations with two different linearizations
are presented by \cite{RPPCalvete2019} for unlimited supply problems with unit-buying customers and uniform budgets. This is extended to a variant of problem in \cite{RPPT} where customers are allowed ties in their ranking of the products and further generalised in \cite{CapicitatedRPP} in which product supply is limited and customers may have different reservation prices for products. In the latter, while an envy-free solution is always possible, the authors consider the implications of solutions with envy in situations where the seller may prefer profit maximisation over customer dissatisfaction.

Multi-unit models, as the name implies, allows for customers to purchase one or more of a collection of heterogeneous items or one or more copies of the same homogeneous good. \cite{ChenNing2008} introduce a notion of metric substitutability in the sense that the same item can be purchased at different locations but will incur a travel cost. For the unit demand model a polynomial time exact algorithm is provided. The extension to customers wishing to purchase up to a maximum amount of items from a location is shown to be APX-hard and a logarithmic approximation algorithm is provided. A variant of the multi-unit model is considered by \cite{Chen2016} in which the seller has items of differing quality and customers require a fixed (sharp) amount of items. Here, the seller wishes to maximise revenue and customers have a known valuation for an item of unit quality. Envy-free and competitive equilibrium solutions are considered, the latter being the requirement where all buyers receive a utility maximising allocation and unallocated items are priced at zero. They develop polynomial time algorithms for specific cases of the two solution types. For homogenous items, \cite{branzei2016envy} consider envy-free solutions when customers have a budget with which to purchase multiple items of the good. They present an FPTAS and exact algorithm for the objective of the seller maximising revenue and a polynomial time algorithm for the objective of maximising social welfare (the overall valuation of items allocated to customers).  \cite{Monaco_RevenueMaxEnvy-freeforhomogeneousResources} introduce the notion of pair envy-freeness, where customer prefer their allocation over other customers' allocations. Customer valuations for bundles of homogeneous items of different sizes can be single-minded, non-decreasing in the number of items or more general where no assumption is made about the valuation relative to the size of the bundle. Hardness and approximation results are presented for the different cases under the two notions of envy-freeness.  \cite{Flammini2019} give a further relaxation of envy-freeness where solutions are  socially envy-free, in which customers prefer their allocation other other customers who are socially connected to them (modelled via some undirected graph of relationships) and provide approximation results of single minded and general valuations

When we consider an underlying network structure to the pricing problem and the motivating real-life problem, \cite{ChenNing2008} for their multi-unit variant, specifically motivate an application of fuel pricing to individual customers choosing amongst all locations based on the price set by the seller and the commute cost to that location. Given that the envy-free pricing problem is hard to approximate, \cite{Guruswami2005} consider the Tollbooth problem as special case. Taking items as edges on a graph (with infinite supply) to model road segments, customers wish to purchase paths (a bundle of items) in the graph and the seller aims to price tolls on these segments to maximise profit. This is shown to be APX-hard. In their LP formulation, \cite{myklebust2016efficient}, utilise customer segments that contain multiples of homogeneous unit-buying customers rather than individual customers. Using the dual of the LP they introduce new and improved heuristics to the pricing problem. We adapt this in our model to allow for segments to be customers with the same route choice but their individual fuel demands may vary. \cite{Fernandes2016} also consider the network pricing problem where tolls are placed on the edges of a graph and customers will pay the sum of tolls from edges used. They give a reduction from the envy-free pricing problem to the basic network pricing problem that preserves the approximation ratio and develop new MIP formulations of the envy-free pricing problem.

\subsection{Contributions}
\begin{itemize}

 \item From conceptual point of view: Our work
    \begin{enumerate}
         \item  shows the importance and impact of explicitly modeling \textit{price match} as a pricing strategy. The conclusion from results in Sections \ref{war_vs_match} and \ref{mnpp_vs_bmnpp} is that when brand sensitivity can be estimated efficiently price matching strategy is very profitable, however, when this is not the case, explicitly modeling it may still very useful but its impact is less severe as price war is most common,
    \item gives a more practical extension of min-pricing, than that is more commonly considered in the previous literature, by incorporating a logit model when a buyer faces choices between two sellers,
    \item illustrates the importance of network structure on pricing decisions, more specifically, our results highlight how locations of sources (fueling stations) within a network influences the pricing decisions and hence the revenue earned by sellers,
    \item illustrates how multi-product pricing models can be applied to network pricing scenarios particularly within fuel pricing contexts.
    \end{enumerate}
   
\item From algorithmic point of view: We
\begin{enumerate}
    \item derive two MIP formulations, to the best of our knowledge, MIP formulations for min-buying are never addressed,
    \item model the case with general demand, as against the unit demand case considered in the existing literature,
    \item propose several heuristics to solve the problem. The main focus of this work from an algorithmic standpoint is towards exploration rather than exploitation. That is, the aim is to set the platform for industrial production level algorithms but not to develop one as such an exercise is out of scope to be covered in one paper. As we point in future directions, the algorithms proposed will form basis for a developing more complex algorithms such as hyper-heuristics within which algorithms from this work constitute as low-level heuristics.
\end{enumerate}
   
\end{itemize}

\subsection{Min-Pricing (MNPP)} \label{sec:MNPP}
In this setting, demand at outlet $f$  from  $e\in \mathcal{N}$ is:
\begin{equation}\label{mnpp-1}
\left\{\begin{array}{ll}
         D_f(p_f,e),  \quad & \mbox{when $p_f\leq p_{f'}, \ f'\in  \mathcal{O}_e$, } \\
     0,\quad & \mbox{otherwise},
       \end{array} \right.
\end{equation}
where 
\begin{align}
    D_f(p_f,e) = \left\{\begin{array}{llr}
        D^{PM}_{ef} = d_e \cdot \beta_e,   \quad & \mbox{when $p_f=c_e$},  & \quad \mbox{(price-match)}\label{mnpp-2}\\
       D^{PW}_{ef} = d_e \cdot \gamma_e,  \quad &  \mbox{when $p_f\leq c^-_e$}, &\quad \mbox{(price-war)}\\
     0, \quad & \mbox{otherwise,} 
       \end{array} \right.
\end{align}
and $0<\beta_e<1$ and $\beta_e < \gamma_e \leq 1$ are constants.
In words, outlet $f$ can satisfy (any) demand at $e$, only if, $f$ is the cheapest priced among $\mathcal{O}_e$. Denoting $A$ as the allocation matrix, that is, $A_{ef}=1$ when $e\in \mathcal{N}$ is assigned to $f\in \mathcal{O}$, and $0$ otherwise. 

Note that $d_e$ is the demand at $e$ which is associated with some route $h\in \mathcal{H}_{(i,j)}$. Hence, in full notation $d_e=d(i,j,h) = \alpha(i,j,h) \cdot d(i,j)$.

\begin{observation}
Given a price ladder $T$, optimal allocation $A$ can be computed by solving a dynamic programme in $O(|\mathcal{O}|\times |M|) $ time.
\end{observation}
\textit{Proof of Observation} \quad First we observe that it is optimal to allocate each demand to first outlet to which it can be allocated to in a price ladder. To see this, let $\mathcal{O}_e=\{T_f, T_{f+r}\}$. Suppose $e$ is allocated to $T_{f+r}$ in one of many possible optimal allocations then since $T_f$ is placed before $T_{f+r}$ in the price ladder and $e$ is allocated to $T_{f+r}$ it must be the case that $T_f$ and $T_{f+r}$ are both given the same price. Therefore, $e$ can be reallocated to $T_f$ to obtain an alternative optimum while keeping the revenue unchanged. 

The above observation gives an easy way of constructing optimal allocation for a given price ladder, however, it only gives an optimal revenue generating allocation but not optimal prices. Let the price index set $M$ be such that prices are increasing, that is, $b(m) \leq b(m+1)$, $1\leq m < |M|$. 

For a given price ladder and the associated allocation, the optimal prices can be computed by the following dynamic programme $DP[\mathcal{O},T,A]$:
\begin{align}
    G(f,m) &= \left (\sum_{e\in \mathcal{N}}  A_{eT(f)} \left ( D_{T(f)}(b(m),e) \cdot b(m) \right) \right ) + \max_{j\leq m} G(T(f-1),j) \quad \forall f>1, m\in M, \\
    G(1,m) &= \sum_{e\in \mathcal{N}} A_{eT(1)} \left (  D_{T(1)}(b(m),e)\cdot b(m) \right).
\end{align}
Note that $G(f,m)$ is the cumulative maximum revenue from the first $f$ outlets when $f^{th}$ in the price ladder, that is, outlet $T(f)$ is priced at $b(m)$.




\subsection{Binary Logit Min pricing (BMNPP)}

In this setting, demand at outlet $f$  from $e\in \mathcal{N}$ is:
\begin{equation}
\left\{\begin{array}{ll}
         D_f(p_f,e),  \quad & \mbox{when $p_f\leq p_{f'}, f'\in  \mathcal{O}_e$, } \\
     0,\quad & \mbox{otherwise},
       \end{array}\right. 
\end{equation}
where 
\begin{align}
    D_f(p_f,e) = \left\{\begin{array}{llr}
         d_e \cdot \frac{\bar{\eta}_{ef}}{1+ \bar{\eta}_{ef}},   \quad & \mbox{when $p_f=c_e$},  & \quad \mbox{(price-match)}\label{bmnpp-2}\\
       d_e \cdot \frac{\eta_{ef}}{1+ \eta_{ef}},  \quad & \bar{c}_e \leq  p_f\leq c_e, &\quad \mbox{(price-war)}\\
     0\quad & \mbox{otherwise}, 
       \end{array}\right. 
\end{align}
where $\bar{c}_e$ is the start of price basket for $e$.\\

In words, as in MNPP, outlet $f$ can satisfy (any) demand at segment $h$ of $e$, only if, $f$ is the cheapest priced among $\mathcal{O}_e$. The key difference between BMNPP and MNPP is as follows. The choice faced by the buyer between the seller and the non-seller outlet is modeled using a binary logit model in both price-match and price-war scenarios. Consider the price war case, that is, if $\bar{c}_e \leq  p_f\leq c_e$ is the price of the seller at outlet $f$ then $\eta_{ef} = exp(\hat{a}_{ef} - \hat{b}_{ef}p_f)$ ($\eta_{ef}$ is usually referred to as preference weight of outlet $f$); where $\hat{a}_{ef} \in (-\infty,\infty)$ and $\hat{b}_{ef} \in [0,\infty]$ are parameters capturing the effect of price on buyers' preference. With this, $\frac{\eta_{ef}}{1+ \eta_{ef}}$ is the choice probability of the seller's outlet $e$.

\begin{remark}
    MNPP is shown to be hard to approximate in previous literature, see \cite{UniformBudget_Briest}. By extension similar complexity results carryover to BMNPP. Similarly, the DP to used solve MNPP under a price ladder can be easily adapted to BMNPP.
\end{remark}
\section{IP formulations}

In this section, two IP formulations are presented for solving MNPP.
In MNPP-IP1, we have following variables: $x_f$ ($x_e$) is price variable for $f\in \mathcal{O}$ ($e\in \mathcal{N}$); $\mu_{ef}$ is a binary variable which takes a value of 1 when $d_e$ is satisfied by $f$. 
\allowdisplaybreaks
\begin{align}
 (\mbox{MNPP-IP1}) \quad \max \quad  \sum_{e\in  \mathcal{N}} (y_e+z_e)    \label{ip1} \\
   s.t. \quad |x_e - x_f| \leq \pi \quad e,f\in (\mathcal{O}\cup \mathcal{N})  \label{ip2}  \\
   \sum_{f \in \mathcal{O}_e} \mu_{e,f} \le 1 \quad e\in \mathcal{N}    \label{ip3}\\
    x_e - x_f \leq \pi(1-\mu_{e,f}) \quad f\in \mathcal{O}_e,e\in\mathcal{N}_f     \\
    x_f - x_e \leq \pi(1-\mu_{e,f}) \quad f\in \mathcal{O}_e,e\in \mathcal{N}_f      \\
     x_{f_1} - x_{f_2} \leq \pi(1-\mu_{e, f_1}) \quad f_1,f_2\in \mathcal{O}_e,e\in \mathcal{N}_f  \label{ip4} \\
      \begin{rcases}
 (c_e - x_e)\leq M w^1_{e} \\
 (x_e-c_e) \leq M w^2_{e} \\
 d_ec_e\beta_e w^3_{e} \geq y_e \leq  x_e d_e \beta_e \\
  w^3_e \le \sum_{f \in \mathcal{O}_e} \mu_{e,f} \\
 w^1_{e}+w^2_{e}+w^3_{e} = 1 \\
 (x_e- c^{-}_e) \leq M w^4_{e} \\
d_e c^{-}_e\gamma_e w^5_{e} \geq z_e \leq  x_e d_e \gamma_e \\
w^4_e \le \sum_{f \in \mathcal{O}_e} \mu_{e,f}\\
 w^4_{e} + w^5_{e} = 1\\
\end{rcases}  \text{ for each $e\in  \mathcal{N}$};    \label{ip5} \\
   \mu_{e, f} \in \{0,1\} \quad e\in \mathcal{N}, f \in \mathcal{O}   \\
  w^i_{e}\in \{0,1\} \quad e\in \mathcal{N}, i=1,2,3,4,5 \\
  x_e, y_e, z_e  \geq 0  \quad e\in \mathcal{N} 
\end{align}
(\ref{ip1}) maximizes revenue through price-match and price-war, (\ref{ip2}) allows for a price differential such that the prices at the seller's outlets can differ by no more than $\pi$, (\ref{ip3})--(\ref{ip4}) ensure that demand $d_e$ is satisfied at a single seller's outlet of least price, (\ref{ip5}) evaluates the revenue from the pricing decisions.

MNPP-IP2 is an extension of \cite{RPPT} to the network min-pricing problem and includes the following variables: $y^{em}_f$, a binary variable, which takes a value of 1 when $d_e$ is satisfied by $f$ at price $m$; $v^m_f$, binary variable, which takes a value 1 if $f\in \mathcal{O}$ is priced at $b(m)$.
\begin{align}
 (\mbox{MNPP-IP2}) \quad \max \quad  \sum_{e\in \mathcal{N}} \sum_{f\in \mathcal{O}_e} \sum_{m\in M} D_f(b(m),e)\cdot y^{em}_f \label{rp1}\\
   s.t. \quad  \sum_{f\in \mathcal{O}_e} \sum_{m\in M} y^{em}_f \leq 1 \quad e\in \mathcal{N} \label{rp2}\\
   \sum_{m\in M} v^{m}_f = 1 \quad f\in \mathcal{O} \label{unique-price}\\
    v^m_f + \sum_{b(\hat{m})<b(m);\hat{m}\in M}\sum_{\hat{f} \neq f; \hat{f}\in \mathcal{O}_e}v^{\hat{m}}_{\hat{f}} \leq y^{em}_f + |\mathcal{O}_e|*(1-y^{em}_f) \quad e\in \mathcal{N},f\in \mathcal{O}_e,  m\in M  \label{rp3}\\
  y^{em}_f \leq v^m_f \quad e\in \mathcal{N},f\in \mathcal{O}_e,  m\in M  \label{rp31}\\
  v^m_f + v^{\hat{m}}_{j} \leq 1 \quad f,j \in \mathcal{O}, m,\hat{m}\in M, b(\hat{m})>b(m)+\pi \label{rp4}\\
   y^{em}_f, v^m_f \in \{0,1\} \quad e\in \mathcal{N}, m \in M, f \in \m{O}_e \label{rp5}
\end{align}
(\ref{rp1}) maximizes revenue, (\ref{rp2}) implies each demand is assigned to at most one outlet, (\ref{rp3}) implies demand $d_e$ purchases at $f$ priced at $b(m)$ only if there no other outlet in $\mathcal{O}_e$ priced at $b(\hat{m})<b(m)$, (\ref{rp4}) is the price differential constraint.



\section{Heuristics - MNPP and BMNPP}
In this section, two classes of constructive heuristics are proposed. The first class is referred to as \textit{selection heuristics} and the second as \textit{insertion heuristics}. Selection heuristics are computationally less intensive compared to the insertion heuristics, but both classes build price ladders to take advantage of polynomial solvability of pricing problem under a price ladder. The benchmark heuristic against which the performance of these heuristics is measured is the pricing mechanism widely used in practice in forecourt fuel pricing due to its simplicity: \textit{Single price heuristic} (SP). Given its name, SP assigns the same price to all outlets and sets the price to be the one that maximises revenue. The proposed algorithms use following two key sub-routines:
 \begin{itemize}
     \item Algorithm \ref{Rev_max_sr} gives revenue from a fixed \textit{allocation}, $A$. 
     \item  Algorithm \ref{allocation_pl_sr} computes an allocation, $A$, of demands to outlets given a price ladder on the outlets.
 \end{itemize}
 Algorithm \ref{single_price_heu} describes the benchmark heuristic which gives the same price to all outlets.

\begin{algorithm}[!htbp]
\caption{Revenue Maximization for a fixed allocation}
\begin{algorithmic}
\State{Input: $\hat{\m{O}} \subseteq \m{O} :=$ set of outlets, $A :=$ allocation matrix, $T :=$ price ladder}
\Function{Revenue}{$A$, $T$, $\hat{\m{O}}$}
\State{$G = DP(A, T, \hat{\m{O}})$}
\State{\Return{$rev=\max_m G(\hat{\mathcal{O}},m)$}}
\EndFunction
\end{algorithmic}
\label{Rev_max_sr}
\end{algorithm}

\begin{algorithm}[!htbp]

\caption{Allocation under price ladder constraint}
\begin{algorithmic}[1]
\Function{Allocate}{$T,nAI$}
\State{$\m{N}' \gets \m{N}$}
\For{$i \in \{1,\dots,nAI\}$}
\State{$f \gets T(i)$}
\For{$e \in \mathcal{N}$}
\If{$f\in \mathcal{O}_e$ and $e \in \m{N}'$}
\State{$A[e][f] \gets 1$}
\State{$\m{N}' \gets \m{N}' \setminus \{e\}$}
\EndIf
\EndFor
\EndFor
\State{\Return{$A$}}
\EndFunction
\end{algorithmic}
\label{allocation_pl_sr}
\end{algorithm}


\begin{algorithm}[!htbp]
\caption{Single Price Heuristics}
\begin{algorithmic}[1]
\State{Input: $\mathcal{N}:= \text{Set of demands}, M:=\text{Set of price indices}$}

\Function{Single-Price}{$\mathcal{N}, M$}
\State{$bestrev \leftarrow -1$; $p \leftarrow 0$}
\For{$m \in M$}
\State{$rev \leftarrow 0$}
\For{$e \in \mathcal{N}$}
\If{$b(m) \le c^-_e$}
\State{$rev \leftarrow rev + \max_{f\in \m{O}_e} b(m) D_f(b(m), e) $}
\EndIf
\EndFor

\If{$bestrev < rev$}
\State{$bestrev \leftarrow rev$; $p \leftarrow b(m)$}
\EndIf
\EndFor
\State{\Return{$p, bestrev$}}
\EndFunction
\end{algorithmic}
\label{single_price_heu}
\end{algorithm}

\subsection{Selection heuristics}\label{selection_heu}

In this class of algorithms, at any iteration of the algorithm a partial price ladder with demand allocated is extended by selecting the next outlet to be inserted in the ladder, hence the name selection heuristics. To emphasize, the $k^{th}$ outlet selected corresponds to the $k^{th}$ outlet in the ladder and will retain its position in all iterations to follow. In other words, a single price ladder is constructed during the course of algorithm. Two algorithms that constitute this class of heuristics are:
\begin{itemize}
    \item Greedy (\textit{greedy}): selects the outlet which produces maximum increase in revenue at each iteration, see Algorithm \ref{greedy_alg},
    \item Outlet ordering (\textit{order}): orders all outlets according to best revenue, see Algorithm \ref{outlet_ordering_alg}. 
\end{itemize}


\begin{algorithm}[!htbp]
\caption{Greedy Heuristic for ordering of Outlets}
\begin{algorithmic}[1]
\State{Input: $\hat{\m{O}} \subseteq \m{O} :=$ set of outlets, $T:= \text{price ladder}$}

\Function{Greedy}{$\hat{\mathcal{O}}$}
\State{$\tilde{\m{N}} \gets \m{N}; \tilde{\m{O}} \gets \hat{\m{O}}; \tilde{\m{O}}_e \gets \hat{\m{O}}_e \ \forall \ e; T \gets (0,\dots, 0); j \gets 1$}
\While{$|\tilde{\m{O}}| > 0$ and $|\tilde{\m{N}}| > 0$}
\State{$\m{N}_f \gets \{e \in \tilde{\m{N}}: f \in \tilde{\m{O}_e}\}$ for $f \in \tilde{\m{O}}$}
\State{$bestrev \gets -1$}
\For{$f \in \tilde{\m{O}}$}
\State{$T(j) \gets f$}
\State{$A_f \gets \text{Allocate}(T, j)$ \Comment{Create allocation based only on the first $j$ outlets}}
\State{$r_f \gets \text{Revenue}(A_f, T, \{T(1),\dots,T(j)\})$}
\If{$r_f < bestrev$}
\State{$f_{\pi} \gets f; bestrev \gets r_f$}
\EndIf
\EndFor
\State{$T(j) \gets f_{\pi}$}
\State{$\tilde{\m{N}} \gets \tilde{\m{N}} \setminus \m{N}_{f_{\pi}}; \tilde{\m{O}} \gets \tilde{\m{O}} \setminus \{f_{\pi}\}$}
\State{$j \gets j + 1$}
\EndWhile
\State{\Return{$T, bestrev$}}
\EndFunction
\end{algorithmic}
\label{greedy_alg}
\end{algorithm}



\begin{algorithm}[!htbp]
\caption{Outlet Ordering Heuristics to order Outlets}
\begin{algorithmic}[1]
\State{Input: $\hat{\m{O}} \subseteq \m{O} :=$ set of outlets, $T:= \text{price ladder}$}

\Function{Order}{$\mathcal{N}, \Tilde{\mathcal{O}}$}
\State{$\tilde{\m{N}} \gets \m{N}; \tilde{\m{O}} \gets \hat{\m{O}}; \tilde{\m{O}}_e \gets \hat{\m{O}}_e \ \forall \ e; T \gets (0,\dots, 0); j \gets 1$}
\While{$|\tilde{\m{O}}| > 0$ and $|\tilde{\m{N}}| > 0$}
\State{$\hat{e} = \arg\min_{e \in \tilde{\m{N}}} c_e$}
\State{$\m{N}_f \gets \{e \in \tilde{\m{N}}: f \in \tilde{\m{O}_e}\}$ for $f \in \tilde{\m{O}}$}
\If{Model $=$ MNPP}
\State{$\mu_f \gets \sum_{e \in \m{N}_f} d_e \gamma_e c_e$ for $f \in \tilde{\m{O}}_{\hat{e}}$}
\EndIf
\If{Model $=$ BMNPP}
\State{$\eta_{e, f} \gets \exp(q_{e, f} - r_{e, f} c_e)$ for $f \in \tilde{\m{O}}_{\hat{e}}, e \in \m{N}_f$}
\State{$\mu_f \gets \sum_{e \in \m{N}_f} d_e \frac{\eta_{e, f}}{\eta_{e, f} + 1} c_e$ for $f \in \tilde{\m{O}}_{\hat{e}}$}
\EndIf
\State{$f_{\pi} \gets \arg\min_{f \in \tilde{\m{O}}_{\hat{e}}} \mu_f$}
\State{$T(j) \gets f_{\pi}$}
\State{$\tilde{\m{O}} \gets \tilde{\m{O}} \setminus \{f_{\pi}\}; \tilde{\m{O}_e} \gets \tilde{\m{O}}_e \setminus \{f_{\pi}\} (e \in \tilde{\m{N}}); \tilde{\m{N}} \gets \tilde{\m{N}} \setminus \m{N}_{f_{\pi}}$}
\EndWhile
\State{\Return{$T$}}
\EndFunction
\end{algorithmic}
\label{outlet_ordering_alg}
\end{algorithm}


\subsection{Insertion Heuristics}
These heuristics have, at their core, an exploration of many price ladders as opposed to the single price ladder in selection heuristics. Recall that, once the price ladder is decided, calculating the revenue maximizing prices is established via DP. These heuristics extend selection heuristics with the selection phase preceding an insertion phase where multiple positions for each outlet are considered at the time of insertion. However, at any iteration previously inserted outlets retain their order in the ladder. It should be noted that insertion heuristics are also constructive in nature.

There are two key steps that make up the skeleton of insertion heuristics: selection and insertion. Selection is done using the selection rules discussed in Section \ref{selection_heu}. An insertion step involves exploring all possible positions that the selected outlet can take in the current price ladder, while preserving the order of outlets previously inserted. The insertion step is formally described in Algorithm \ref{insert}.

A total of five insertion heuristics are devised:
\begin{enumerate}
    \item Full Insertion (\textit{FI}): In the selection step, we consider all outlets still to be inserted, and the best revenue option is picked. 
    \item Greedy heuristic as selection rule followed by insertion (\textit{greedyI}).
    \item Outlet-ordering as selection rule followed by insertion (\textit{orderI}).
    \item An IP1 based selection rule followed by insertion (\textit{IP1I}): the LP-relaxation of IP1 is used to order the outlets based on the price variables which is then used as selection rule.
    \item An IP2 based selection rule followed by insertion (\textit{IP2I}): the LP-relaxation of IP2 is used to order the outlets based on the price variables which is then used as selection rule.
\end{enumerate}

\begin{algorithm}[ph]
\caption{Insert routine}
\label{insert}
\begin{algorithmic}
\State{Input: $\hat{\m{O}} \subseteq \m{O} :=$ set of outlets, $T:= \text{price ladder}$}
\State{Initialize: $T[i]=-1$ for all $i<|\mathcal{O}|$; $\hat{\m{O}}=\{\}$, $nAI=0$}
\Function{Insert}{$f$}
\State{$bestpos(f)=-1$, $bestrev(f)=-1$}
\State{$temp = T$} 
\State{$\hat{\m{O}} = \hat{\m{O}} \cup f$; $nAI \gets nAI + 1$}
\For{$j\in [0,nAI]$}
\State{$tT[i]=-1$ for all $i$}
\For{$i<nAI$}
\If{$i>j$}
\State{$tT[i]=temp[i-1]$}
\Else
\State{$tT[i]=T[i]$}
\EndIf
\State{$tT[j]=f$}
\EndFor
\State{$A(j)$ = Allocate$(tT,nAI)$}
\State{$r(j)=\mbox{Revenue}(A(j))$}
\If{$bestrev(f) < r(j)$}
\State{$bestrev(f) = r(j)$; $bestpos(f)=j$}
\EndIf
\EndFor
\State{$tT[i]=-1$ for all $i$}
\For{$i<nAI$}
\If{$i>bestpos(f)$}
\State{$tT[i]=temp[i-1]$}
\Else
\State{$tT[i]=T[i]$}
\EndIf
\State{$tT[bestpos(f)]=f$}
\EndFor
\State{$T=tT$}
\EndFunction
\For{$f\in \mathcal{O}$ and $f\notin \hat{\m{O}}$ }
\State{Insert$(f)$}
\EndFor
\State{Output: $T$, $bestrev(|\mathcal{O}|)$}
\end{algorithmic}
\end{algorithm}






Briefly, the time complexity of the heuristics is as follows: selection heuristics have a worst-case running time of $O\left(\frac{|\mathcal{O}|\times (|\mathcal{O}|+1)}{2} \times |M|\right)$, whereas, insertion heuristics take $O\left(\left(\frac{|\mathcal{O}|\times (|\mathcal{O}|+1)}{2}\right)^2 \times |M|\right)$.
\section{Numerical experiments}

The aim of this section is to conduct an extensive analysis of the effectiveness of exact and heuristic algorithms using simulated and realistically generated instances. All experiments have been coded in python.

\subsection{Metrics}

To enable comparison between heuristics, apart from the usual gap to optimality metric,
\begin{equation}
   OptGap =  \frac{r(H) - r(OPT)}{r(OPT)} \times 100,
\end{equation}
the following metric will also be employed: 
\begin{equation}
   GS (H) =  \frac{r(H) - r(SP)}{r(SP)} \times 100.
\end{equation}
Gain-over-SP ($GS$)  measures percentage gain of any heuristic that takes network aspect in to account compared to ignoring it. This comparison with SP is motivated by the following reasons:
\begin{itemize}
    \item Theoretically, SP already achieves the best possible approximation, under certain conditions,
    \item Arguably, SP is the simplest possible pricing mechanism, computational cost-wise, which implies an increase in computational complexity can only be justified with the proportional increase in revenue.
\end{itemize}

\subsection{Experimental Design}

For each demand model MNPP and BMNPP, we test our algorithms in 450 instances. We generate these instances as follows. We consider $|\m{O}| \in \{5, 10, 15\}$ and $|\m{N}| \in \{15, 30, 50\}$. The density of the bipartite graph $G$ is defined as $P \in \{0.9, 0.75, 0.5, 0.25, 0.1\}$, where $P = \frac{\lvert L \rvert}{\lvert \m{O} \rvert \times \lvert \m{N} \rvert}.$ For each instance, the set of edges $L$ is chosen randomly, such that the resulting graph has density $P$. In addition, for all MNPP instances, we used $\beta_e = 0.5$ and $\gamma_e = 1$ for all $e \in \m{N}$. This generates 45 bipartite graphs. For each graph, we create 10 instances by sampling the values of $c_e$ uniformly from the interval $[0, 25]$ and $d_e$ uniformly from the interval $[50, 150]$. This generates 450 instances for each demand model. Each algorithm was given a maximum run time of 4 hours. For the logit pricing model, for each $e \in \m{N}$ and $f \in \m{O}$, $\hat{a}_{ef}$ and $\bar{a}_{ef}$ were sampled uniformly from the interval $[200,400]$ and $\hat{b}_{ef}$ and $\bar{b}_{ef}$ were sampled uniformly from $[0, 20]$. Note that we select $\hat{b}_{ef}$ to be much smaller than $\hat{a}_{ef}$ in order to avoid $\eta_{ef}$ becoming infinitesimally small when $p_f$ is not large.
 
\subsection{Results}

\subsubsection{MNPP}\label{sec:MNPP_results}

\begin{figure}[htbp!]
    \centering
    \includegraphics[width=0.40\textwidth]{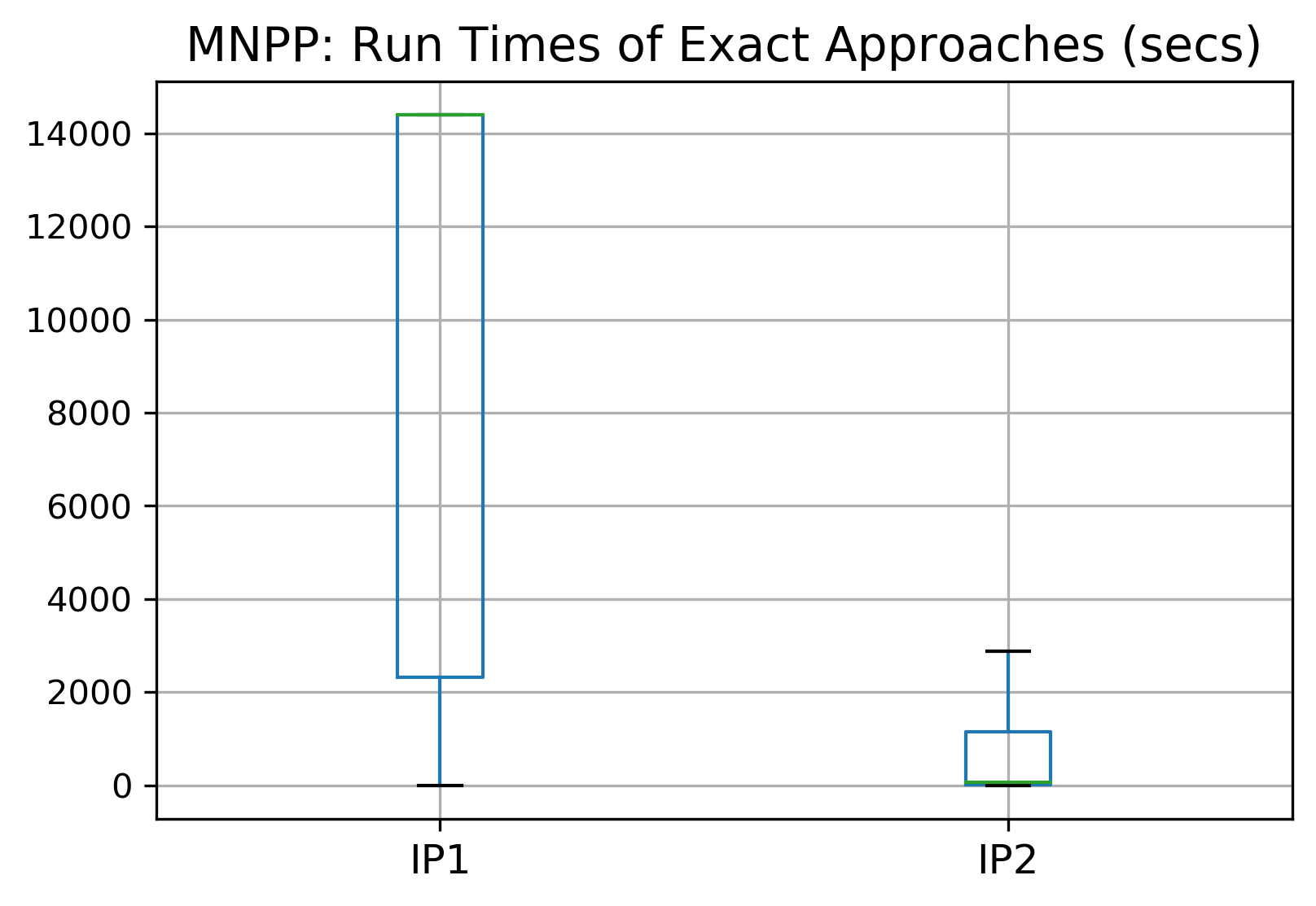}
    \hspace{1cm}
    \includegraphics[width=0.40\textwidth]{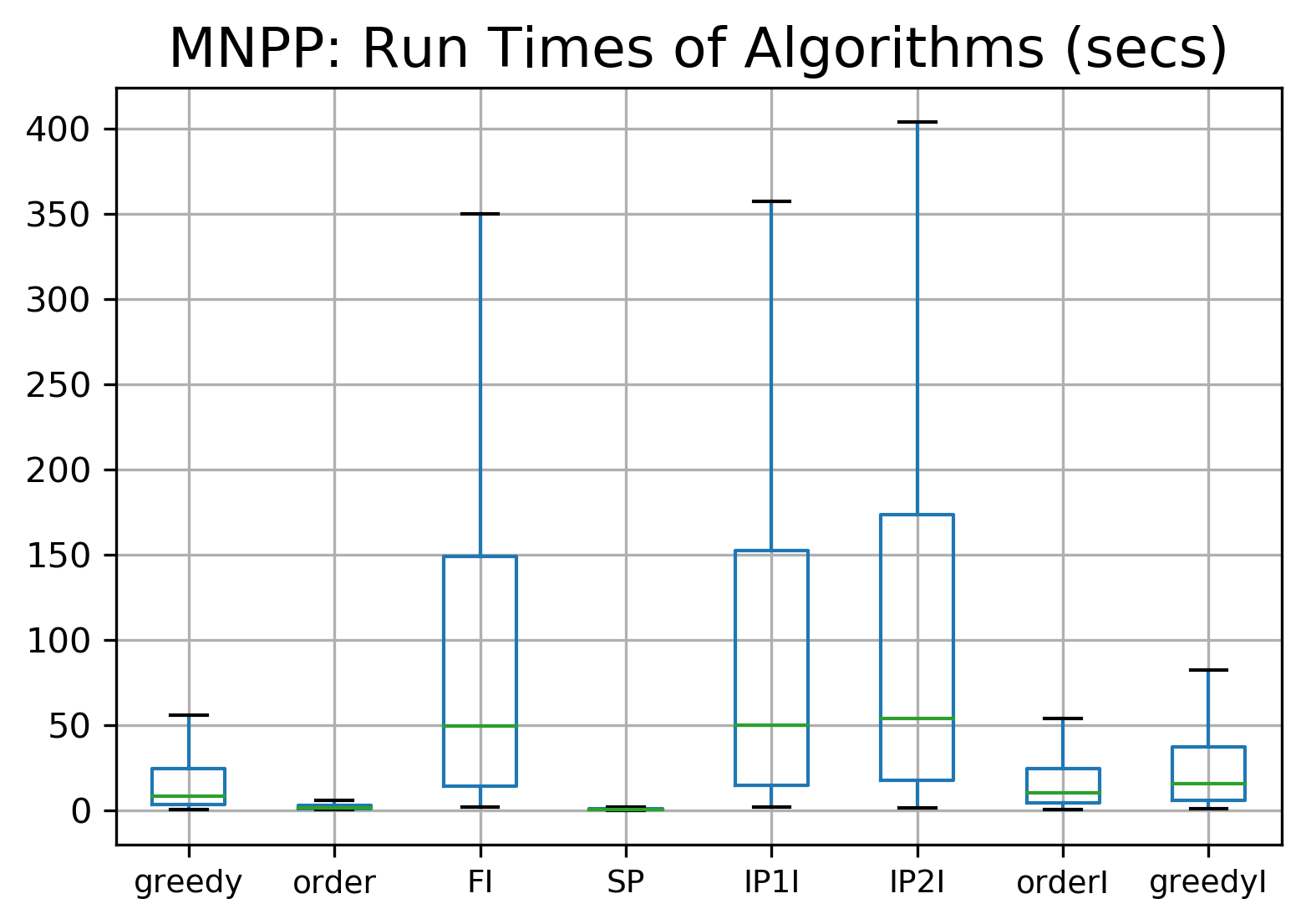}
    \caption{Boxplots summarising times taken by the exact and heuristic algorithms for MNPP.}
    \label{fig:min_times}
\end{figure}

We first present the run times of each algorithm under the MNPP demand model in Figure~\ref{fig:min_times}. Firstly, the left plot in Figure~\ref{fig:min_times} shows that IP2 is significantly faster than IP1. IP1 timed in 296 instances out of 450, whereas IP2 only timed out in 43. Including timeouts, IP1's average run time was 10233 seconds (approximately 2 hours 50 minutes), whereas IP2's was only 2258.4 seconds (approximately 37 minutes). The second plot in Figure~\ref{fig:min_times} shows the run times of our heuristics. Standout results are that FI, IP1I and IP2I are significantly slower than the other algorithms. FI took an average of 1 minute and 57 seconds to finish running. While this is a substantial reduction compared the exact approaches, FI was still slow compared to the other heuristics. SP and order were the fastest heuristics, taking an average of 0.59 and 7.35 seconds respectively. In addition, orderI and greedy offered similar time savings over IP1 and IP2, and both took approximately 20 seconds on average.

\begin{figure}[htbp!]
    \centering
    \includegraphics[width=0.40\textwidth]{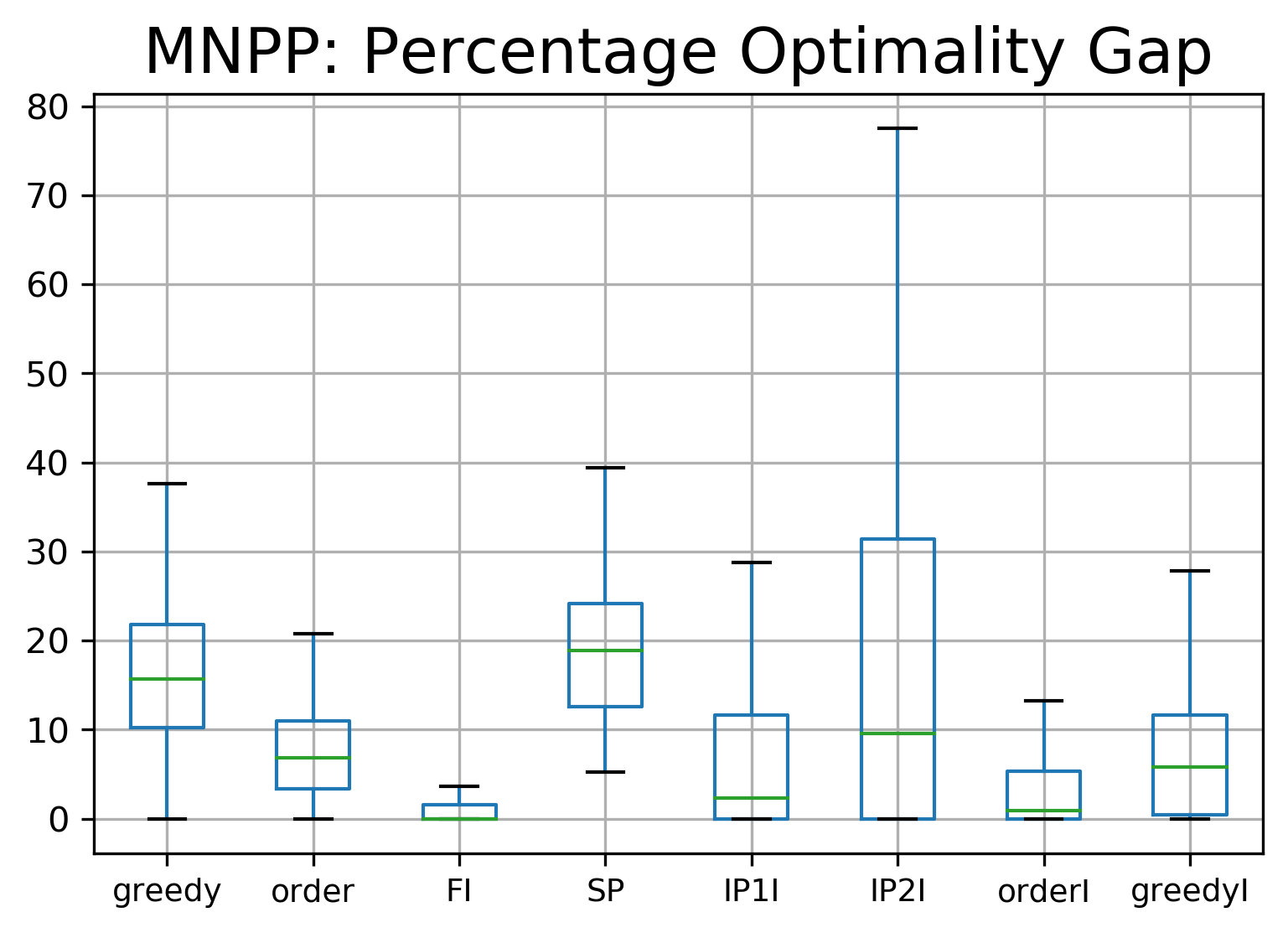}
    \hspace{1cm}
    \includegraphics[width=0.40\textwidth]{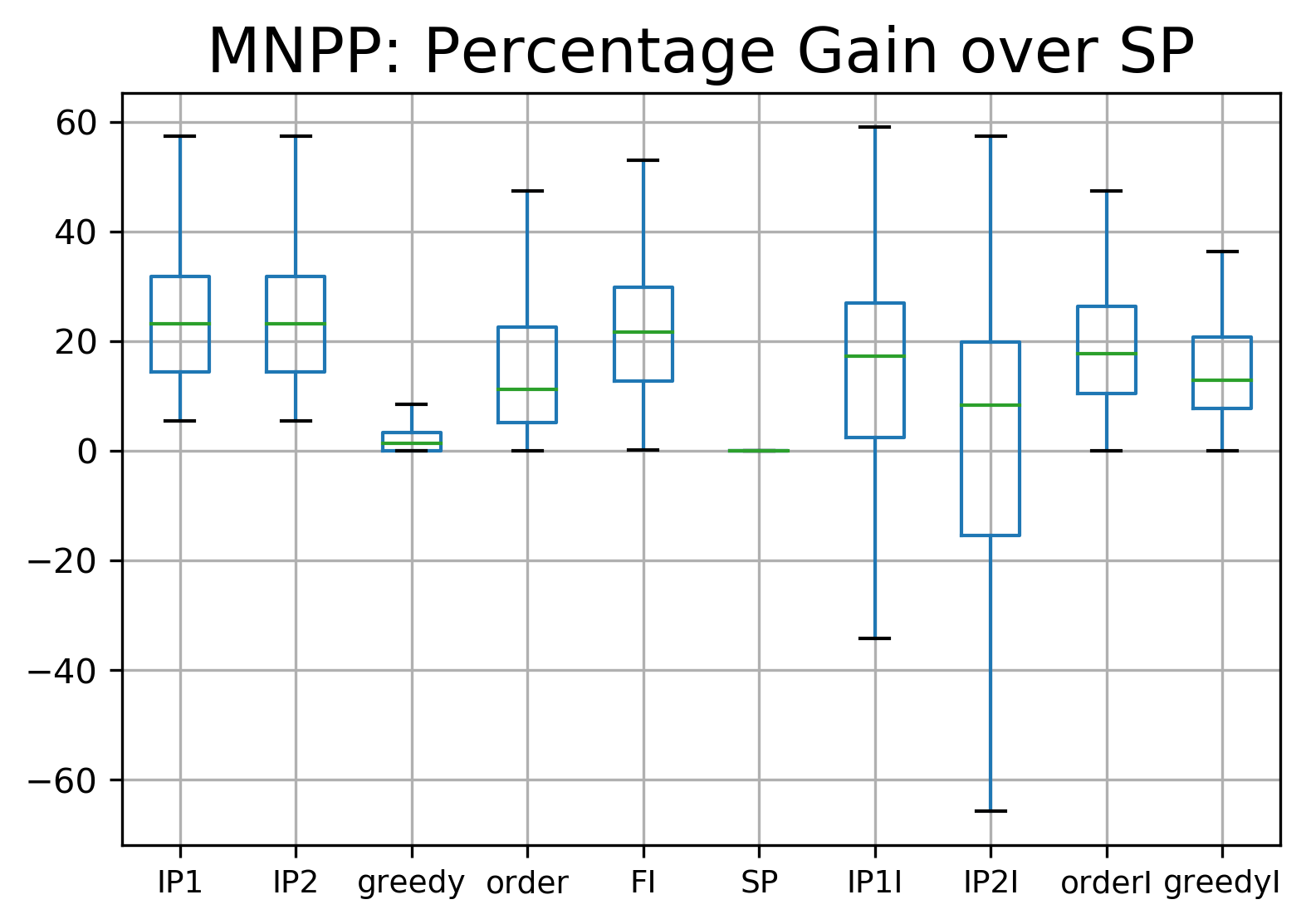}
    \caption{Boxplots summarising percentage optimality gaps and gains over SP for MNPP.}
    \label{fig:min_pog}
\end{figure}

Figure~\ref{fig:min_pog} shows the percentage optimality gaps (left) and gain over SP (right). The optimality gaps indicate that FI was the best performing heuristic for the MNPP demand model. Its mean percentage optimality gap was 1.16\% and it was optimal in 40\% of instances. We also find that order and orderI performed well, achieving average percentage optimality gaps of 7.86\% and 3.6\% respectively. Since these algorithms were significantly faster than FI, it might be more practical to use them over FI. IP1I also performed quite well, with an average percentage optimality gap of 8.07\%. However, since it was slower than order and orderI, there is no reason to use it here. While SP was the fastest heuristic, it had an average percentage optimality gap of 19.5\%, which is much higher than order and orderI. The right-hand plot in Figure~\ref{fig:min_pog} shows the gain over SP of each algorithm. As might be expected, IP1 and IP2 gained the most. It is interesting that IP1I and IP2I also gained as much as IP1 and IP2 in some instances, but they were also outperformed by SP by as much as 60\% in some instances. 

\begin{figure}[htbp!]
    \centering
    \includegraphics[width=0.40\textwidth]{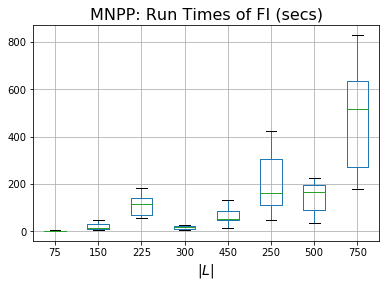}
    \hspace{1cm}
    \includegraphics[width=0.40\textwidth]{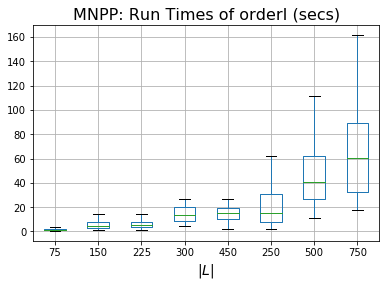}
    \caption{Boxplots summarising times taken by the FI and orderI by $|L|$ for MNPP.}
    \label{fig:min_times_edges}
\end{figure}

Figure~\ref{fig:min_pog} suggests that FI was the best performing heuristic, but Figure~\ref{fig:min_times} suggests that it was significantly slower than order and orderI. Since these algorithms perform similarly to FI, the time difference may make them the preferred heuristics. In order to explain why FI was slow, Figure~\ref{fig:min_times_edges} shows the run times of FI and orderI by $|L|$. The lefthand plot shows that FI took longer than 6 minutes for $|L| = 250$ and over 12 minutes for $|L| =  750$. On the other hand, orderI typically did not take longer than 2 minutes and 45 seconds, even for $|L| = 750$. These results indicate that FI scales poorly with $|L|$ compared to orderI. This is likely because, in every iteration, FI runs the insert algorithm for every single remaining outlet in order to select which one to insert. On the other hand, orderI runs the simple ordering heuristic once in order 
to select the next outlet to insert. This requires significantly less computation and enumeration. 

\begin{figure}[htbp!]
    \centering
    \includegraphics[width=0.8\textwidth]{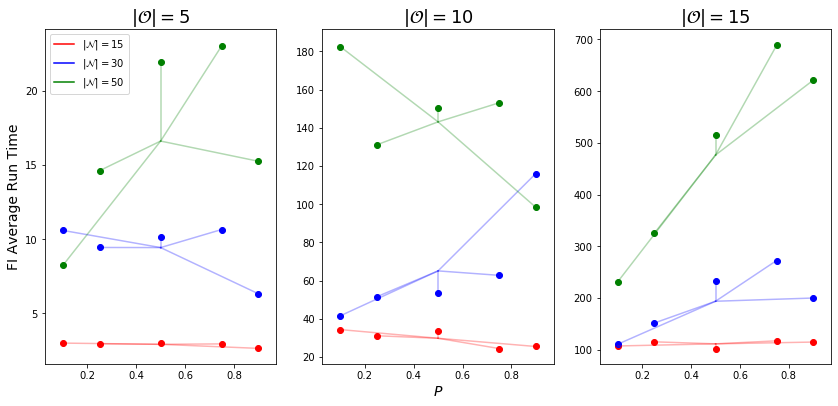}
    \caption{Run times of FI by $|\m{O}|$, $|\m{N}|$ and $P$}
    \label{fig:fancy_FI_tt}
\end{figure}

Figure~\ref{fig:min_times_edges} shows that FI scales poorly with $|L|$ compared to orderI. However, $|L|$ is influenced by 3 other model parameters: $|\m{O}|$, $|\m{N}|$ and $P$. In order to understand how FI and orderI scale with each of these parameters individually and which has the largest effect, we present Figure~\ref{fig:fancy_FI_tt}. The first thing to note is that the effect of $P$ increased significantly with $|\m{N}|$. The value of $P$ had little effect on run times for $|\m{N}| = 15$ for all values of $|\m{O}|$. However, as $|\m{N}|$ increases we begin to see significant more variation in run times over $P$ for each value of $|\m{O}|$. Since density had no effect on times unless $|\m{N}| > 15$, this indicates that density only caused large run times when combined with large $|\m{N}|$. Secondly, Figure~\ref{fig:fancy_FI_tt} indicates that $|\m{O}|$ had a slightly larger effect on FI's run times than $|\m{N}|$. While the effects were similar, the variation in times over $|\m{O}|$ for each value of $|\m{N}|$ was larger, with run times for $|\m{N}| =50$ increasing from around 15 seconds at $|\m{O}| =5$ to over 450 seconds at $|\m{O}| = 15$. Conversely, for $|\m{O}| =15$, times only varied between approximately 100 seconds and 450 seconds. This result is intuitive, since larger $|\m{O}|$ means more outlets and hence more runs of the insert heuristic in each iteration of FI.  

\begin{figure}[htbp!]
    \centering
    \includegraphics[width=0.8\textwidth]{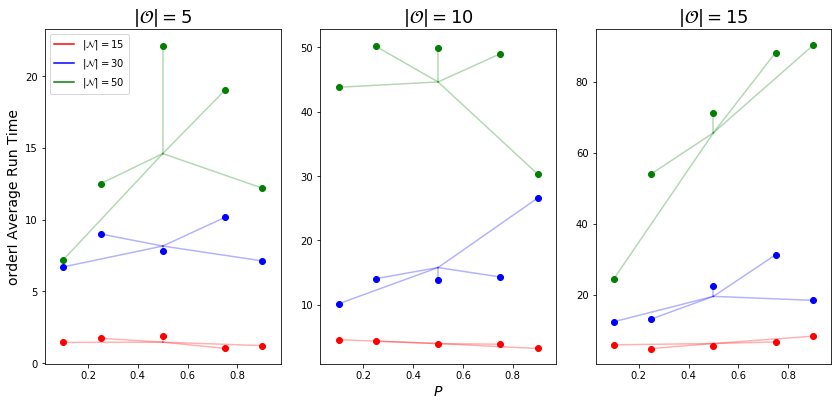}
    \caption{Run times of orderI by $|\m{O}|$, $|\m{N}|$ and $P$}
    \label{fig:fancy_orderI_tt}
\end{figure}

Figure~\ref{fig:fancy_orderI_tt} shows orderI's run times. Firstly, we see that the same conclusion about the effect of $P$ holds here, in that it had no real effect unless $|\m{N}| > 15$. In addition, the effect of $|\m{O}|$ on orderI's run times was much weaker than on FI's run times. Even for the largest values of $|\m{N}|$, orderI's run times only increased from around 15 seconds to around 60 seconds when $|\m{O}|$ is increased from 5 to 15. This is because, where FI runs one insert algorithm for every outlet in $\m{O}$ in every iteration, orderI simply runs order once in each iteration. Hence, it scales much better with $|\m{O}|$. In fact, $|\m{N}|$ had a slightly larger effect on orderI's run times than $|\m{O}|$, with its average time for $|\m{N}| = 50$ being 41 seconds and its average run time for $|\m{O}| = 15$ being only 30 seconds. 

\begin{figure}[htbp!]
    \centering
    \includegraphics[width=0.8\textwidth]{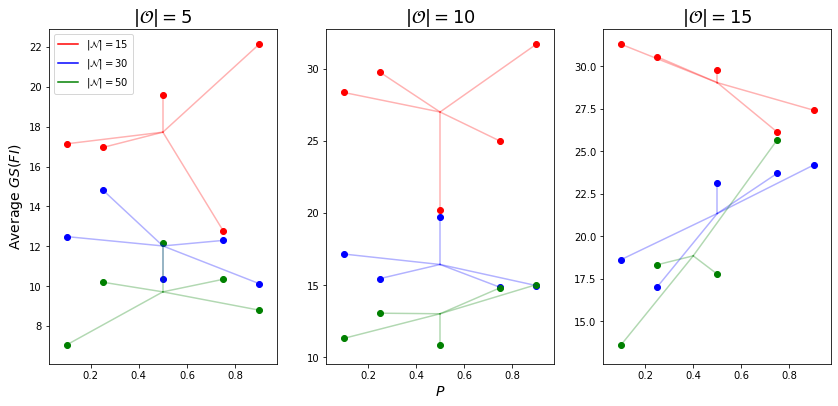}
    \caption{Gain over SP of FI by $|\m{O}|$, $|\m{N}|$ and $P$}
    \label{fig:fancy_FI_gs}
\end{figure}

\begin{figure}[htbp!]
    \centering
    \includegraphics[width=0.8\textwidth]{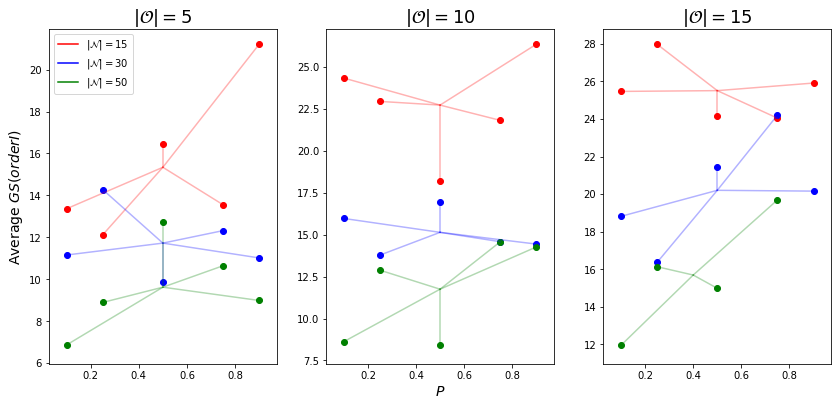}
    \caption{Percentage gains over SP of orderI by $|\m{O}|$, $|\m{N}|$ and $P$}
    \label{fig:fancy_orderI_gs}
\end{figure}

We did not find that $|\m{O}|$, $|\m{N}|$ or $P$ had any significant effect on the percentage optimality gaps of FI or orderI. However, they did have an effect on the amount that each algorithm gained over SP, i.e.\ $GS(H)$ for $H \in \{\text{FI}, \text{orderI}\}$. Figures~\ref{fig:fancy_FI_gs} and~\ref{fig:fancy_orderI_gs} show $GS(\text{FI})$ and $GS(\text{orderI})$ in more detail. As might be expected, the larger the value of $|\m{O}|$, the more revenue the two algorithms gained over SP. Both algorithms gained around 5-15\% additional revenue when $|\m{O}|$ was increased from 5 to 15. This is intuitive due to the following reason. Suppose we introduce a new outlet and the demand nodes remain the same. If the new outlet's demand nodes have high competitor prices, then SP will keep its prices the same, since setting a high price at this new outlet would mean increasing all prices and hence losing demand at other nodes. On the other hand, if the new outlet's demand nodes have low prices, then SP will have to reduce all of its prices to gain demand at the new outlet. This means losing revenue at other outlets. Conversely, orderI and FI can set a high price at the new outlet in the first case and a low price at the new outlet in the second case without affecting the revenue at other nodes. Therefore, they naturally gain more from the increase in outlets. 

Interestingly, however, large values of $|\m{N}|$ had the opposite effect on the algorithms' gains over SP. In general, both algorithms gained more over SP when $|\m{N}|$ was smaller. Both algorithms gained an additional 10-15\% revenue over SP when $|\m{N}|$ was 15 compared to when it was 50. This result is due to the following reason. Suppose we have a fixed set of outlets and we add a new demand node. If SP's prices are already cheaper than this demand's competitors, then its prices will stay the same. However, if not, then SP will need to reduce all of its prices in order to gain revenue from the new demand node. This can result in a loss of revenue at other outlets. As we increase the number of demand nodes, SP's prices will therefore typically decrease. FI and orderI, however, may only need to reduce one price in order to gain revenue from this node, without affecting revenue at the others.  

In summary, FI can be used for problems on small graphs or when optimality is the most important factor for the decision maker. However, when $|\m{O}|$ is large, FI begins to run very slowly. We also conclude that orderI provides solutions that are almost as good as FI's, but in a shorter time. This is due to the fact that its run times are much less affected by $|\m{O}|$ than FI's. We note that FI and orderI's performance are not particularly affected by the structure of the graph, but that incorporating this information results in 25-30\% gains over SP when there are a large number of outlets and few demands. We also found that orderI and FI gain more over SP when there are fewer outlets. Finally, if speed is of the most importance, one can use order, which ran in approximately a third of the time that orderI ran in, but gave solutions that were slightly worse. 


\subsubsection{BMNPP}\label{sec:BMNPP_results}

\begin{figure}[htbp!]
    \centering
    \includegraphics[width=0.40\textwidth]{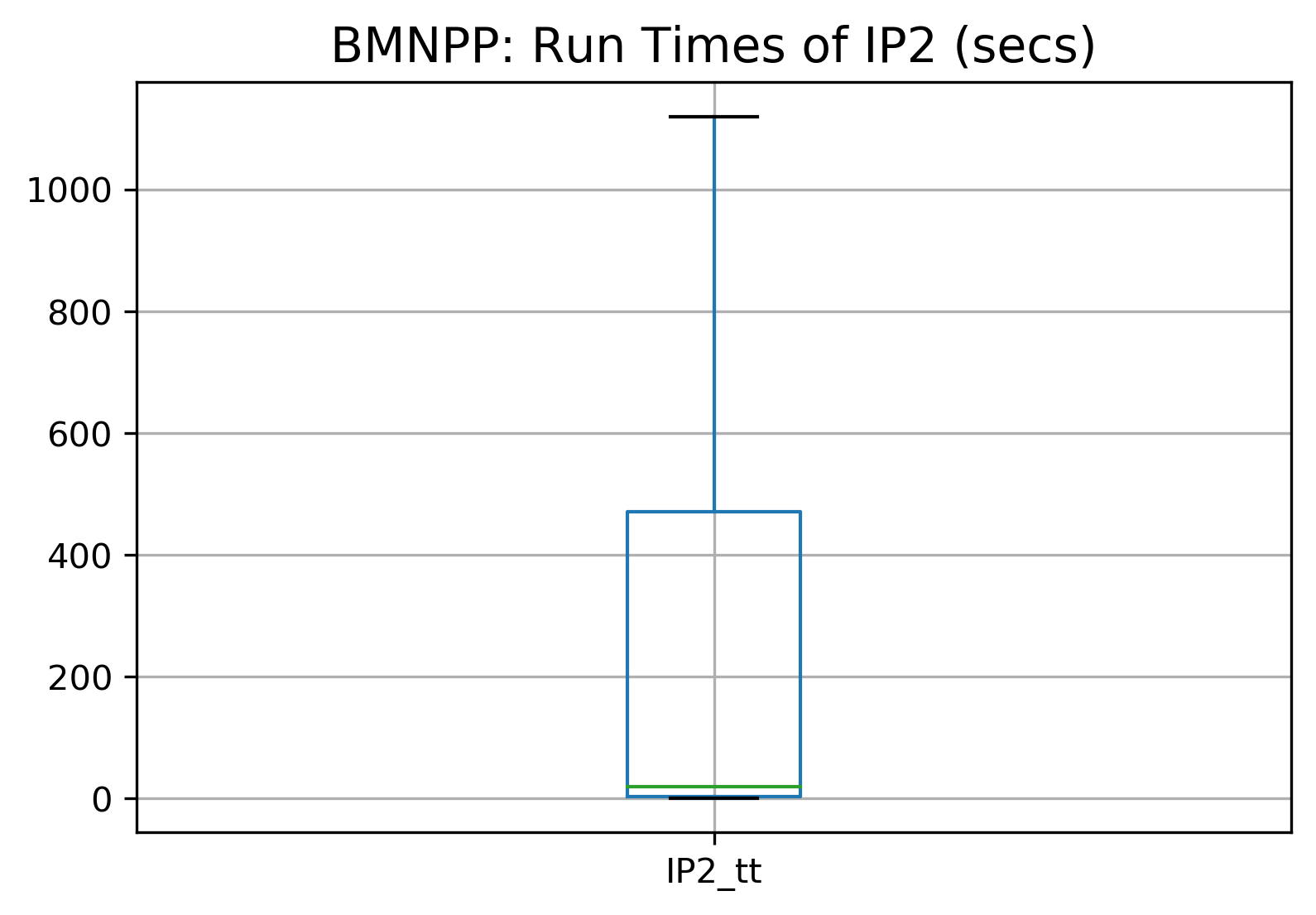}
    \hspace{1cm}
    \includegraphics[width=0.40\textwidth]{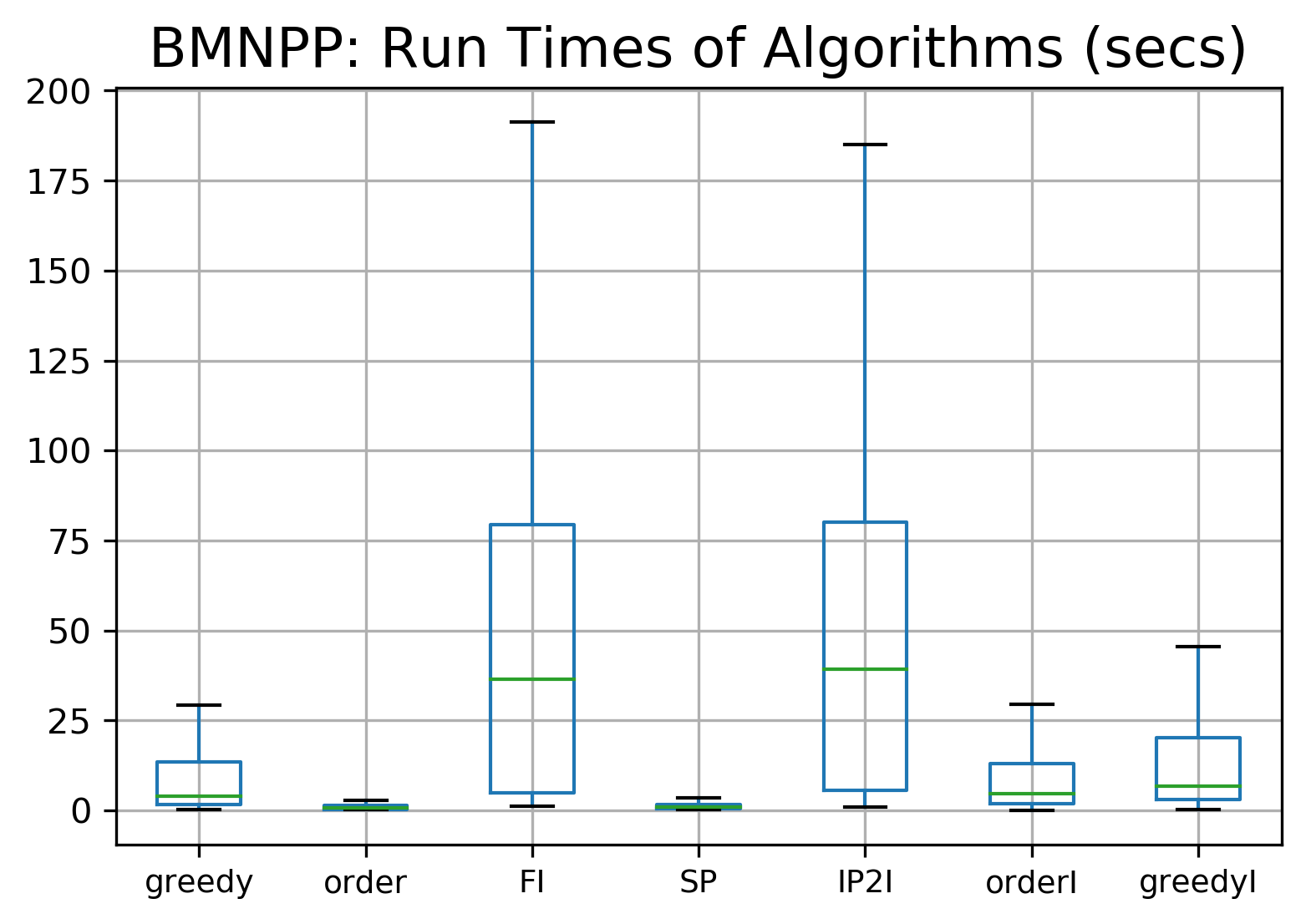}
    \caption{Boxplots summarising times taken by the exact and heuristic algorithms for BMNPP.}
    \label{fig:logit_times}
\end{figure}

We now present the results for BMNPP. Note that we did not run IP1 for BMNPP as the results for MNPP indicate that there is no reason to use it over IP2. Figure~\ref{fig:logit_times} shows the run times of all algorithms. The left-hand figure shows that IP2 was faster for BMNPP than MNPP, which is likely due to the fact that the exponential demand model means that the demands become zero for high prices. This allows us to eliminate high prices as potential choices, reducing the number of decision variables. For these experiments, IP2 only timed out in 25 instances as opposed for 43 for MNPP. However, it is still the slowest algorithm by far, taking an average of 24 minutes and 1 second to finish running. The right-hand plot in Figure~\ref{fig:logit_times} shows the run times of our heuristics. As before, FI and IP2I were the slowest heuristics, while SP and order were the fastest. FI took an average of 1 minute and 4 seconds to run for BMNPP while SP and order took 1.33 and 1.08 seconds respectively.  In addition, orderI was faster for BMNPP, taking only 9.75 seconds on average.

\begin{figure}[htbp!]
    \centering
    \includegraphics[width=0.40\textwidth]{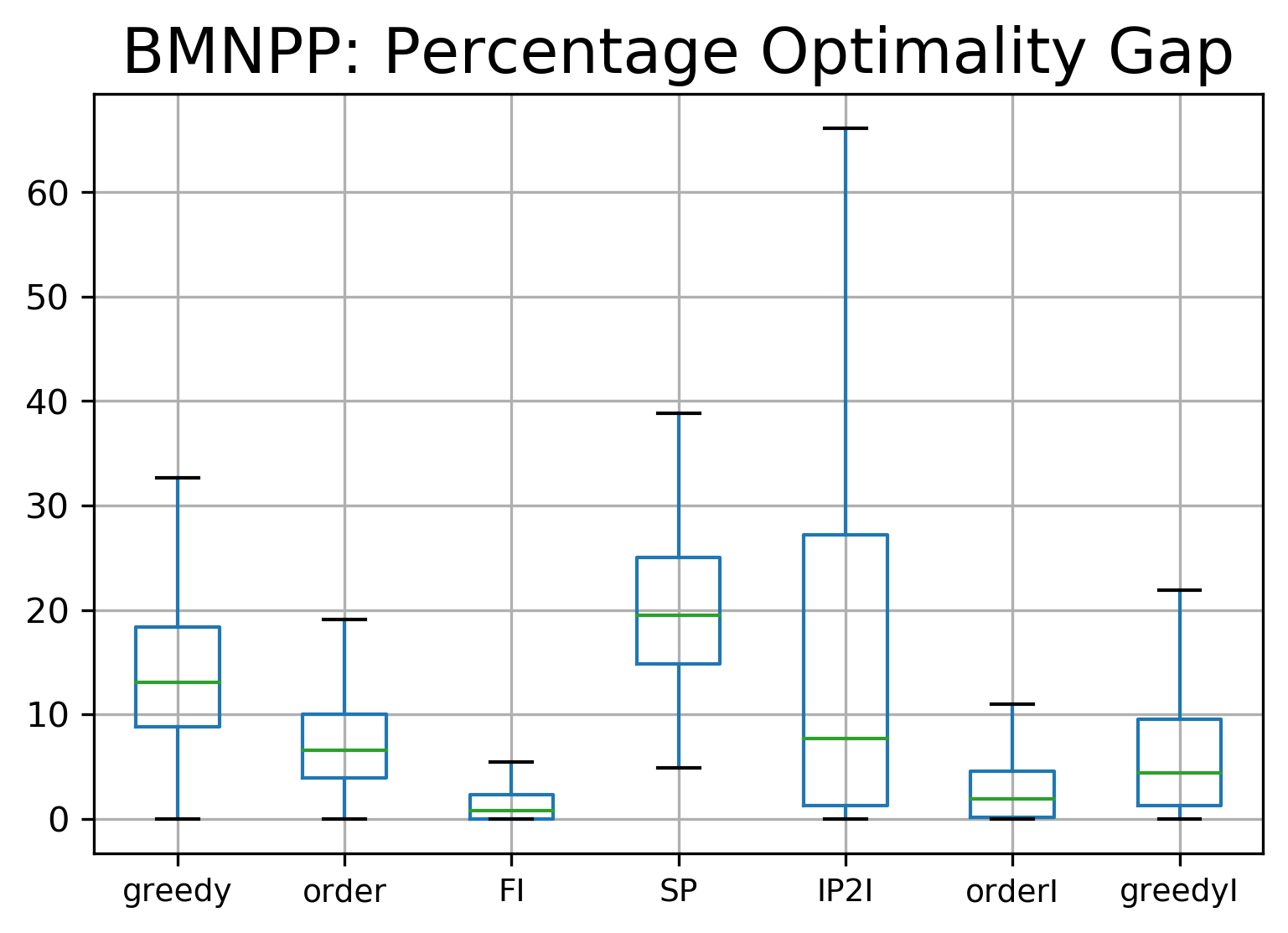}
    \hspace{1cm}
    \includegraphics[width=0.40\textwidth]{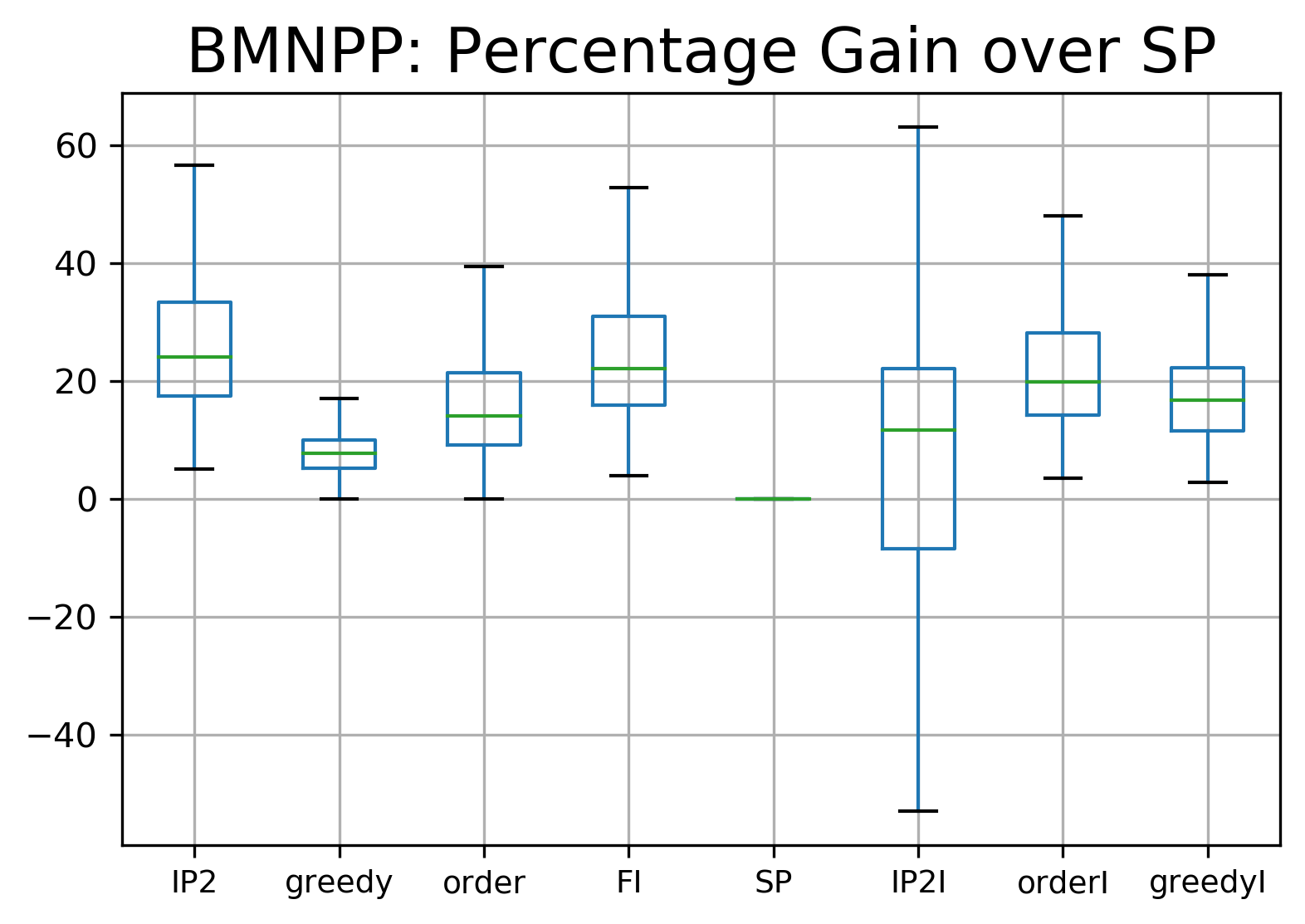}
    \caption{Boxplots summarising percentage optimality gaps and gains over SP for BMNPP.}
    \label{fig:logit_pog}
\end{figure}

The heuristics also performed similarly for BMNPP as they did for MNPP. In particular, the left-hand plot in Figure~\ref{fig:logit_pog} shows that FI was again the best algorithm, closely followed by the much-faster orderI. FI's average percentage optimality gap was 1.47\% while orderI's was 3.12\%. We also find that order had an average percentage optimality gap of 7.4\%, indicating again that it is a fast alternative that provides only slightly worse solutions. In addition, greedyI performs slightly better than order, with an average percentage optimality gap of 6.1\%. However, its average run time was 17.8 seconds as opposed to the better-performing orderI's average time of 9.75 seconds. Therefore, there is no reason to use greedyI over orderI.  As before, SP was the fastest heuristic, but its high average percentage optimality gap of 20\% makes it impractical. As shown in the right-hand plot, IP2, order and orderI gained the most by considering network structure. They typically gained around 20\% revenue over SP as a result of this. We also see that the only algorithm that did not always gain revenue from considering the network structure was IP2I. This is likely due to the effects of relaxing the integality requirements on the decision variables.
\\

\begin{figure}[htbp!]
    \centering
    \includegraphics[width=0.40\textwidth]{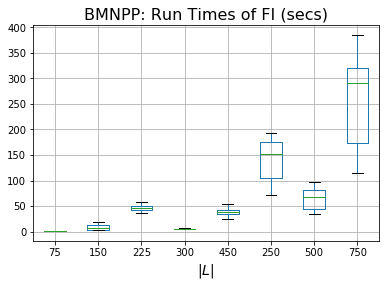}
    \hspace{1cm}
    \includegraphics[width=0.40\textwidth]{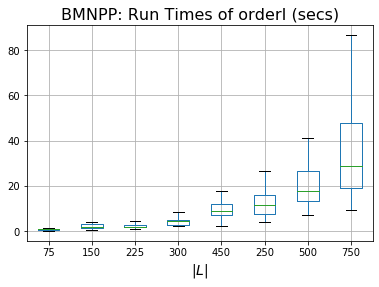}
    \caption{Boxplots summarising times taken by the FI and orderI by $|L|$ for BMNPP.}
    \label{fig:logit_times_edges}
\end{figure}

In summary, the best performing heuristics were FI and orderI, but orderI was much faster than FI. Figure~\ref{fig:logit_times_edges} shows the run times of FI and orderI by $|L|$. This plot indicates that, again, FI scales poorly with $|L|$, its times reaching above 6 minutes for the most connected graphs. While orderI still becomes slower as $|L|$ increases, it still did not take significantly more than 1 minute and 20 seconds in these instances.

\begin{figure}[htbp!]
    \centering
    \includegraphics[width=0.8\textwidth]{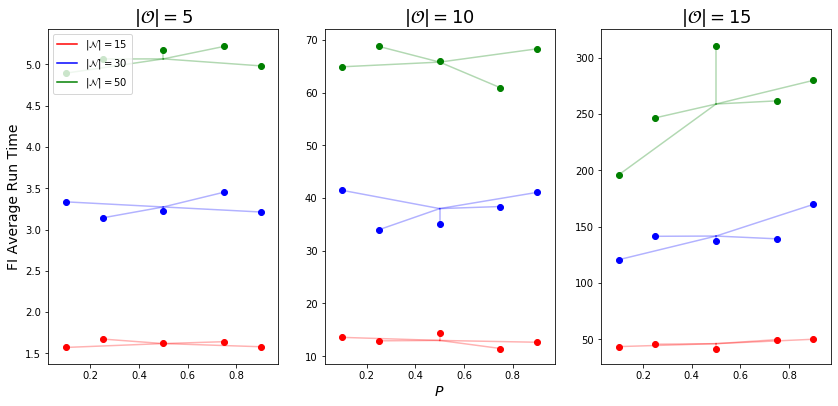}
    \caption{Run times of FI by $|\m{O}|$, $|\m{N}|$ and $P$}
    \label{fig:logit_fancy_FI_tt}
\end{figure}

\begin{figure}[htbp!]
    \centering
    \includegraphics[width=0.8\textwidth]{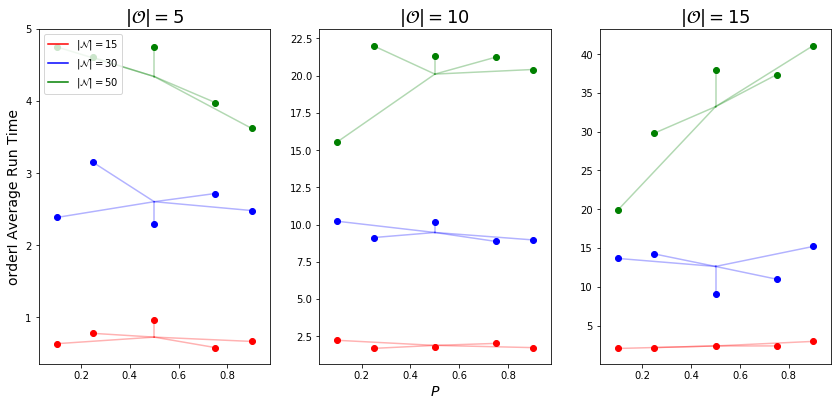}
    \caption{Run times of orderI by $|\m{O}|$, $|\m{N}|$ and $P$}
    \label{fig:logit_fancy_orderI_tt}
\end{figure}

 In order to explain the effect of $|L|$ on run times and solution quality, we now present scatter plots of these run times and GS broken down by $|\m{N}|$, $|\m{O}|$ and $P$. Firstly, Figures~\ref{fig:logit_fancy_FI_tt} and~\ref{fig:logit_fancy_orderI_tt} summarise FI and orderI's run times. Similarly to MNPP, we observe that $P$ had no real effect on the run times of either algorithm when $|\m{N}| = 15$. Its effect for large $|\m{N}|$ is less apparent for BMNPP than for MNPP, but both algorithms were the most affected by $P$ when $|\m{N}| = 50$. In addition, the increase in FI and orderI's run times from increasing $|\m{N}|$ is more clear here, with larger relative differences between the collections of points for each $|\m{O}|$. Figure~\ref{fig:logit_fancy_FI_tt} suggests that, again, $|\m{O}|$ had the strongest effect on FI's run times, with the times for the largest $|\m{N}|$ increasing drastically with each increase in $|\m{O}|$. On the other hand, Figure~\ref{fig:logit_fancy_orderI_tt} suggests that the effects of $|\m{N}|$ and $|\m{O}|$ on orderI's run times were quite similar. As before, orderI's run times were much less affected by $|\m{O}|$ than FI's.

\begin{figure}[htbp!]
    \centering
    \includegraphics[width=0.8\textwidth]{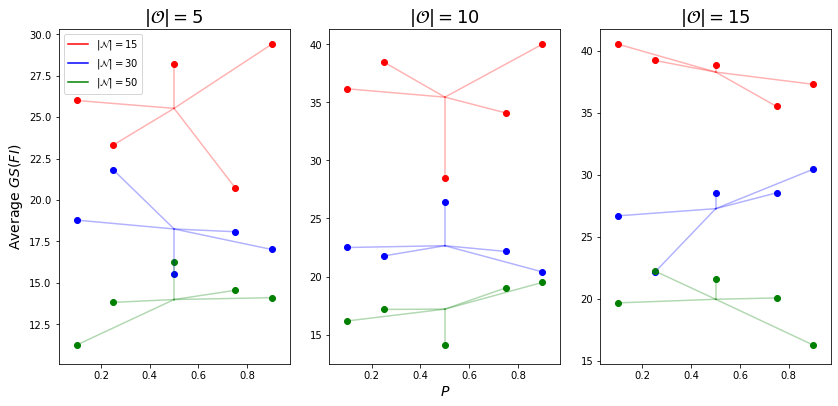}
    \caption{Percentage gains over SP of FI by $|\m{O}|$, $|\m{N}|$ and $P$}
    \label{fig:logit_fancy_FI_gs}
\end{figure}

\begin{figure}[htbp!]
    \centering
    \includegraphics[width=0.8\textwidth]{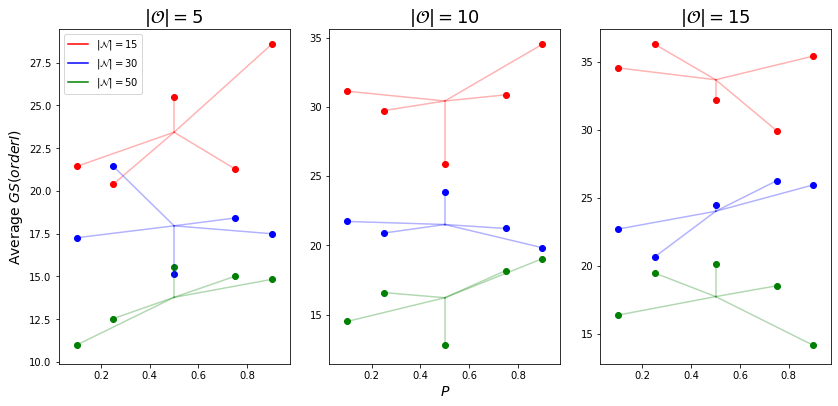}
    \caption{Percentage gains over SP of orderI by $|\m{O}|$, $|\m{N}|$ and $P$}
    \label{fig:logit_fancy_orderI_gs}
\end{figure}

Figures~\ref{fig:logit_fancy_FI_gs} and~\ref{fig:logit_fancy_orderI_gs} show $GS(H)$ for $H \in \{\text{FI}, \text{orderI}\}$.  Both plots indicate again that FI and orderI gain much more over SP when $|\m{N}|$ is small compared to when it is large. In fact, they gain around additional 17-25\% revenue over SP when $|\m{N}| = 15$ compared to when $|\m{N}| = 50$. This is an even bigger difference than for MNPP, which is likely due to the fact that demand under the BMNPP model is also dependent on price. In addition, both algorithms again gain more over SP when $|\m{O}|$ is larger. The biggest difference here occurs when $|\m{N}| = 15$, where FI gains an additional 10-17\% revenue over SP when $|\m{O}| = 15$ compared to when $|\m{O}| = 5$. A similar effect is found for orderI, although orderI only gains slightly less additional revenue when $|\m{O}|$ increases.

Similarly to for MNPP, we conclude that FI can be used for small graphs or in instances where optimality is of the most importance. However, it will begin to solve slowly when $|\m{O}|$ becomes large. Otherwise, orderI can be used to gain strong solutions in a shorter time for all sizes of graph. We also conclude that both FI and orderI gain the most over SP when $|\m{N}|$ is small and $|\m{O}|$ is large. Finally, order can be used to make significant reductions in run times but lose slightly more revenue.

\subsubsection{Price Match vs.\ Price War}\label{war_vs_match}

We now assess the amount of demand achieved under each model by price matching and by price warring. In particular, suppose we obtain two sets of prices, $p^{MNPP}_f$ and $p^{BMNPP}_f$ for $f \in \m{O}$. For pricing model $x \in \{MNPP, BMNPP\}$, define the set of demands that $f \in \m{O}$ price matched and price warred under model $x$ as follows:
\begin{equation}
    \m{N}^{x}_f(PM) = \left\{e \in \m{N}_f: p^x_f = c_e\right\}, \quad \m{N}^{x}_f(PW) = \left\{e \in \m{N}_f: p^x_f \le c^-_e\right\}.
\end{equation}
Then, we calculate the demand obtained by model $x \in \{MNPP, BMNPP\}$ resulting from price matching and price warring via:
\begin{equation}
    d^x(PM) = \sum_{f \in \m{O}}\sum_{e \in \m{N}^{x}_f(PM)} D_f(p^x_f, e), \quad d^x(PW) = \sum_{f \in \m{O}}\sum_{e \in \m{N}^{x}_f(PW)} D_f(p^x_f, e)
\end{equation}
The total demand obtained over all outlets under model $x$ is therefore $d^x = d^x(PM) + d^x(PW)$ and the percentage of demand obtained from price matching and warring respectively can be calculated using $d^x_{\%}(PM) = 100\times\frac{d^x(PM)}{d^x}$ and $d^x_{\%}(PW) = 100\times\frac{d^x(PW)}{d^x}$ respectively. We can similarly define the percentages of total revenue obtained from price matching and price warring under each model $x$ as $r^x_{\%}(PM), r^x_{\%}(PW)$.

Since price matching is a new addition to the network pricing problem, we will now analyse how often MNPP and BMNPP made use of price matching and how much revenue they obtained from it. Table~\ref{tab:MNPP_pm} summarises MNPP's price matching behaviour.  Note that the demand described here is the demand obtained by MNPP's solution under MNPP's demand model. In particular, we notice that over all instances MNPP obtained only 8.85\% of its demand from price matching. The average number of price matches by outlets under MNPP was only 6.09. This small number of price matching resulted in only 10.53\% of MNPP's revenue being obtained from price matching. We also see that MNPP began to obtain more of its demand from price matching as either $|\m{O}|$ or $|\m{N}|$ increased, with the highest average value of 15.41\% occurring when both of these parameters were largest.

\begin{table}[htbp!]
    \centering
    \begin{tabular}{llrrr}
\toprule
        &    &  Average PM Count &  Average $d^{MNPP}_{\%}(PM)$ &  Average $r^{MNPP}_{\%}(PM)$ \\
$|\mathcal{N}|$ & $|\mathcal{O}|$ &                   &                              &                              \\
\midrule
15 & 5 &              1.06 &                         3.82 &                         4.22 \\
        & 10 &              2.84 &                         9.08 &                        11.23 \\
        & 15 &              3.62 &                        10.45 &                        11.97 \\
30 & 5 &              2.94 &                         5.01 &                         5.70 \\
        & 10 &              6.26 &                        10.12 &                        12.20 \\
        & 15 &              9.12 &                        13.36 &                        15.98 \\
50 & 5 &              4.73 &                         5.15 &                         5.61 \\
        & 10 &             10.37 &                         9.30 &                        11.29 \\
        & 15 &             17.85 &                        15.41 &                        19.32 \\
        \hline 
Overall &    &              6.09 &                         8.85 &                        10.53 \\
\bottomrule
\end{tabular}

    \caption{Summary of MNPP's price matching}
    \label{tab:MNPP_pm}
\end{table}

\begin{figure}[htbp!]
    \centering
    \includegraphics[width=0.6\textwidth]{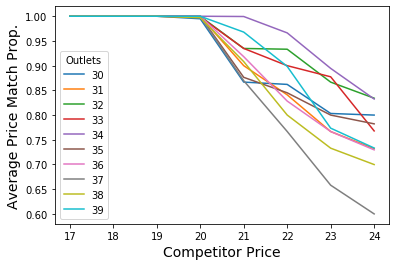}
    \caption{Average price match proportions for BMNPP}
    \label{fig:mvl_props}
\end{figure}

\begin{table}[htbp!]
    \centering
    \begin{tabular}{llrrr}
\toprule
        &    &  Average PM Count &  Average $d^{BMNPP}_{\%}(PM)$ &  Average $r^{BMNPP}_{\%}(PM)$ \\
$|\mathcal{N}|$ & $|\mathcal{O}|$ &                   &                               &                               \\
\midrule
15 & 5 &              5.16 &                         27.39 &                         30.81 \\
        & 10 &              8.40 &                         39.32 &                         46.67 \\
        & 15 &             10.22 &                         46.90 &                         54.28 \\
30 & 5 &              7.67 &                         16.74 &                         18.62 \\
        & 10 &             13.78 &                         28.32 &                         33.38 \\
        & 15 &             17.38 &                         33.90 &                         41.10 \\
50 & 5 &              9.57 &                         12.61 &                         13.85 \\
        & 10 &             18.37 &                         20.43 &                         24.28 \\
        & 15 &             23.81 &                         25.26 &                         30.89 \\
        \hline 
Overall &    &             12.13 &                         28.10 &                         32.85 \\
\bottomrule
\end{tabular}

    \caption{Summary of BMNPP's price matching}
    \label{tab:BMNPP_pm}
\end{table}
Table~\ref{tab:BMNPP_pm} shows that BMNPP utilised price matching significantly more than MNPP. On average, BMNPP used price matches twice as often as MNPP did. This is due to the fact that BMNPP accounts for the fact that price matching can result in more than 50\% of demand from a given demand node and that price warring does not always result in 100\% of this demand. To see this, we studied some instances where BMNPP price matched at least 20 times while MNPP price matched less than 10 times. Figure~\ref{fig:mvl_props} shows the average price match demand proportions over all demand nodes for each outlet. We see that, in this instance, BMNPP's pricing model suggests that price matching will always net us at least 60\% of demand on average. This proportion only decreases from 1 when the competitor price becomes higher than 18. On the other hand, MNPP's model would have assumed that many of these price matches would only result in 50\% of demand, hence undervaluing price matching significantly.   

The effect of price matching is reflected in the average demand and revenue values. Under MNPP, price matching an average of 6.09 times per instance was assumed to result in only 8.85\% of demand being obtained from price matching due the simple 50\% rule. Under BMNPP, however, price matching typically results in more demand than this, which is reflected by the fact that price matching 12.13 times per instance resulted in 28.10\% of total demand being obtained from price matching. Accounting for the fact that price matching can result in more demand from a node than 50\% meant that BMNPP was able to set higher prices than MNPP at a large number of outlets, and this ultimately resulted in more revenue.

In general, we can conclude that the BMNPP pricing model promotes price matching much more than MNPP. Typically, MNPP will inspire us to refrain from price matching as it suggests that this will result in a large loss of demand. However, using BMNPP allows us to account for the fact that more than $50\%$ of customers can buy our fuel even if we do price match with competitors. The fact that BMNPP's price match demand proportions vary among outlets also allows us to distinguish between outlets that need to price war and those that don't as opposed to assuming that they all do, resulting in more opportunities for higher prices.

\subsubsection{MNPP vs.\ BMNPP}\label{mnpp_vs_bmnpp}

In this section, we compare the solutions resulting from the two different pricing models. In particular, we assess the optimality gap of the solution resulting from MNPP under the BMNPP model. This allows us to understand the revenue loss from assuming an MNPP model over the more complex BMNPP model. In particular, we describe:
\begin{equation}
    OptGap(MNPP) = 100 \times \frac{r^*_{BMNPP} - r_{BMNPP}(MNPP)}{r^*_{BMNPP}},
\end{equation}
where $r^*_{BMNPP}$ is the optimal revenue under the BMNPP pricing model. We summarise the values of $OptGap(MNPP)$ for the 421 instances in which IP2 did not time out under MNPP or BMNPP in Table~\ref{tab:min_pog}. Table~\ref{tab:min_pog} shows that the average value of $OptGap(MNPP)$ was 5.92\%, and from the 75\% quantile we see that its value was not often larger than 6.77\%. On average, ignoring graph structure does not result in a significant loss of revenue relatively speaking. This is likely due to MNPP's tendency to price war. Since MNPP sets its prices lower than BMNPP, it gains more demand than BMNPP but receives lower prices for fuel. Therefore, using MNPP results in a more aggressive pricing strategy that sets low prices for in order to gain high demand.

\begin{table}[htbp!]
    \centering
    \begin{tabular}{lr}
\toprule
{} &  $OptGap(MNPP)$ \\
\midrule
Count &       421 \\
Mean  &         5.92\% \\
Std   &         3.38\%\\
Min   &         1.24\% \\
25\% Quantile   &         3.89\% \\
50\%  Quantile &         5.60\%\\
75\% Quantile  &         6.77\% \\
Max  &        28.28\% \\
\bottomrule
\end{tabular}

    \caption{Summary of $OptGap(MNPP)$ values.}
    \label{tab:min_pog}
\end{table}

\begin{figure}[htbp!]
    \centering
    \includegraphics[width=0.8\textwidth]{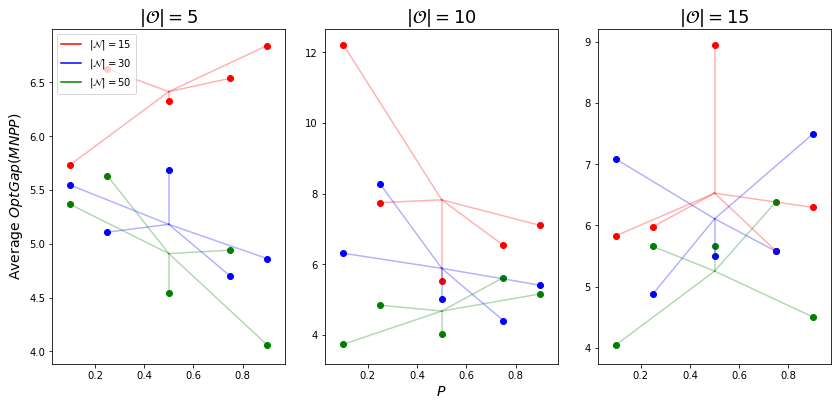}
    \caption{$OptGap(MNPP)$ broken down by $|\m{O}|$, $|\m{N}|$ and $P$}
    \label{fig:fancy_min_pog}
\end{figure}

Despite the average $OptGap(MNPP)$ being only 5.92\%, on occasion this value was as high as 28.28\%. The value of $OptGap(MNPP)$ was higher than 20\% in 8 instances, with varying values of $|\m{O}|$ and $|\m{N}|$. Therefore, it is clear that in some cases the more aggressive pricing of MNPP can lose larger amounts of revenue. In order to analyse $OptGap(MNPP)$ further, we present a plot similar to those in Sections~\ref{sec:MNPP_results} and~\ref{sec:BMNPP_results} that help us observe the effect of $P$, $|\m{N}|$ and $|\m{O}|$ on $OptGap(MNPP)$. Figure~\ref{fig:fancy_min_pog} firstly indicates that $OptGap(MNPP)$ has no clear relationship with $P$. For the some combinations of values of $|\m{O}|$ and $|\m{N}|$, $OptGap(MNPP)$ seems to reduce as $P$ increases. Conversely, for other combinations, $OptGap(MNPP)$ increases with $P$ or does not have a clear relationship with $P$. One result that is clear, however, is that $OptGap(MNPP)$ was the highest when $|\m{N}|$ was lowest, for every value of $|\m{O}|$. This indicates that assuming that the demand function is independent of price results in a higher loss of revenue when there are fewer demand nodes. 

Our analysis suggests that this result is due to the differences in price matching between the two models. We find that as $|\m{N}|$ increased, BMNPP price matched a lower and lower proportion of the demand nodes whereas MNPP price matched slightly more demand nodes. When $|\m{N}| = 15$, we found that BMNPP price matched 52\% of demand nodes whereas MNPP price matched only 17\% of demand nodes. For comparison, when $|\m{N}| = 50$ we found that BMNPP and MNPP price matched 32\% and 20\% of demand nodes respectively. This result indicates that price matching has more benefit when the number of demand nodes is lower, which occurs due to the following reason. When a new demand node is introduced, if the new node has expensive competitors then outlets will not price match, since they would then lose all demand from nodes that were previously price matched. On the other hand, if the new node has cheap competitors, some outlets may reduce their prices in order to price match the new node. When there are a large number of demand nodes, this may result in significant losses in revenue from the old nodes for a slight increase in revenue from one only node. It is more likely that, in this case, outlets will leave their prices the same and receive demand from the new node via price warring.  Therefore, as the number of demand nodes increases, BMNPP begins to behave more like MNPP, in that it price matches a lower proportion of demand nodes and begins to prefer to price war the new demand nodes. 
\section{Conclusion and Future work}

In this work a pricing problem over networks is studied inspired by practical pricing situations such as forecourt retail fuel pricing. The model proposed captures several practical aspects often encountered such as \textit{price match}, \textit{non-uniform demand} and \textit{brand sensitivity} which is ignored in the literature. 

From the algorithmic point of view, this work focuses on exploration of algorithmic possibilities to efficiently solve MNPP and BMNPP. Specifically, two integer programming formulations and two classes of heuristics are analyzed. In terms of future work, embedding these proposed heuristics in an advanced hyper-heuristic framework will produce an algorithm which can be efficiently deployed on industrial scale instances. However, this is not in the scope of this study in order to keep it within practical writing limits. Such an embedding will require discussion on various acceptance and stopping criterion along with extended experimentation for tuning the algorithm for optimal performance. However, we found the best optimality gaps of proposed algorithms are relatively small on average, hence, employing a heavy algorithmic mechanisms such hyper-heuristics may not yield huge benefits. A second promising direction is to reduce running times by reducing calls to DP or algorithmic improvements to speed-up DP. 

Conceptually, a promising extension, especially with increasing electric vehicle uptake, would be to accommodate electric vehicle charging within the model. This will need modeling charging times and their impact on route selection probabilities and capacity constraints at charging stations. The proposed model can accommodate these settings by introducing a time cost associated with charging and waiting times. 

While designed to operate at the operational pricing level of decision-making, our models can be embedded into strategic level network design decisions, such as determining the location and number of outlets to operate within the network.



\normalsize
\bibliographystyle{apalike}
\bibliography{main}

\begin{thebibliography}{}

\bibitem[Aggarwal et~al., 2004]{aggarwal2004algorithms}
Aggarwal, G., Feder, T., Motwani, R., and Zhu, A. (2004).
\newblock Algorithms for multi-product pricing.
\newblock In {\em International Colloquium on Automata, Languages, and Programming}, pages 72--83.

\bibitem[Br{\^a}nzei et~al., 2016]{branzei2016envy}
Br{\^a}nzei, S., Filos-Ratsikas, A., Miltersen, P.~B., and Zeng, Y. (2016).
\newblock Envy-free pricing in multi-unit markets.
\newblock {\em arXiv preprint arXiv:1602.08719}.

\bibitem[Briest, 2008]{UniformBudget_Briest}
Briest, P. (2008).
\newblock Uniform budgets and the envy-free pricing problem.
\newblock In {\em International Colloquium on Automata, Languages, and Programming}, pages 808--819.

\bibitem[Briest and Krysta, 2007]{BuyingCheapisexpensive_Briest}
Briest, P. and Krysta, P. (2007).
\newblock Buying cheap is expensive: Hardness of non-parametric multi-product pricing.
\newblock In {\em Proceedings of the Eighteenth Annual {ACM-SIAM} Symposium on Discrete Algorithms}, pages 716--725.

\bibitem[Calvete et~al., 2019]{RPPCalvete2019}
Calvete, H.~I., Dom{\'\i}nguez, C., Gal{\'e}, C., Labb{\'e}, M., and Marin, A. (2019).
\newblock The rank pricing problem: Models and branch-and-cut algorithms.
\newblock {\em Computers \& Operations Research}, 105:12--31.

\bibitem[Chen et~al., 2016]{Chen2016}
Chen, N., Deng, X., Goldberg, P.~W., and Zhang, J. (2016).
\newblock On revenue maximization with sharp multi-unit demands.
\newblock {\em Journal of Combinatorial Optimization}, 31(3):1174--1205.

\bibitem[Chen et~al., 2008]{ChenNing2008}
Chen, N., Ghosh, A., and Vassilvitskii, S. (2008).
\newblock Optimal envy-free pricing with metric substitutability.
\newblock In {\em Proceedings of the 9th ACM Conference on Electronic Commerce}, pages 60--69.

\bibitem[Dom{\'\i}nguez et~al., 2021]{RPPT}
Dom{\'\i}nguez, C., Labb{\'e}, M., and Mar{\'\i}n, A. (2021).
\newblock The rank pricing problem with ties.
\newblock {\em European Journal of Operational Research}, 294(2):492--506.

\bibitem[Dom{\'\i}nguez et~al., 2022]{CapicitatedRPP}
Dom{\'\i}nguez, C., Labb{\'e}, M., and Mar{\'\i}n, A. (2022).
\newblock Mixed-integer formulations for the capacitated rank pricing problem with envy.
\newblock {\em Computers \& Operations Research}, 140:105664.

\bibitem[Fernandes et~al., 2016]{Fernandes2016}
Fernandes, C.~G., Ferreira, C.~E., Franco, A.~J., and Schouery, R.~C. (2016).
\newblock The envy-free pricing problem, unit-demand markets and connections with the network pricing problem.
\newblock {\em Discrete Optimization}, 22:141--161.

\bibitem[Flammini et~al., 2019]{Flammini2019}
Flammini, M., Mauro, M., and Tonelli, M. (2019).
\newblock On social envy-freeness in multi-unit markets.
\newblock {\em Artificial Intelligence}, 269:1--26.

\bibitem[Guruswami et~al., 2005]{Guruswami2005}
Guruswami, V., Hartline, J.~D., Karlin, A.~R., Kempe, D., Kenyon, C., and McSherry, F. (2005).
\newblock On profit-maximizing envy-free pricing.
\newblock In {\em Proceedings of the Sixteenth Annual {ACM-SIAM} Symposium on Discrete Algorithms}, pages 1164--1173.

\bibitem[Monaco et~al., 2015]{Monaco_RevenueMaxEnvy-freeforhomogeneousResources}
Monaco, G., Sankowski, P., and Zhang, Q. (2015).
\newblock Revenue maximization envy-free pricing for homogeneous resources.
\newblock In {\em Twenty-Fourth International Joint Conference on Artificial Intelligence}, pages 90--96.

\bibitem[Myklebust et~al., 2016]{myklebust2016efficient}
Myklebust, T.~G., Sharpe, M., and Tun{\c{c}}el, L. (2016).
\newblock Efficient heuristic algorithms for maximum utility product pricing problems.
\newblock {\em Computers \& Operations Research}, 69:25--39.

\bibitem[Rusmevichientong et~al., 2006]{PaatRusmevichientong}
Rusmevichientong, P., Van~Roy, B., and Glynn, P.~W. (2006).
\newblock A nonparametric approach to multiproduct pricing.
\newblock {\em Operations Research}, 54(1):82--98.

\end{thebibliography}

\newpage

\end{document}